\def\ZZ         {{\bf Z}}
\def\RR         {{\bf R}}
\def\CC         {{\bf C}}
\def\QQ         {{\bf Q}}
\def\PP         {{\bf P}}

\def\A          {{\cal A}}
\def\F          {{\cal F}}
\def\J          {{\cal J}}
\def\L          {{\cal L}}
\def\O          {{\cal O}}  
\def\M          {{\cal M}}
\def\P          {{\cal P}}
\def\S          {{\cal S}}
\def\U          {{\cal U}}

\def\dim        {{\rm dim}}

\def\Supp       {{\rm Supp}}
\def\Spec       {{\rm Spec}}
\def\Coker      {{\rm Coker}}
\def\Ker        {{\rm Ker}}
\def\Char       {{\rm Char}}
\def\Sing       {{\rm Sing}}

\documentstyle[twoside,12pt]{article}
\newtheorem{prop}{Proposition}[section]
\newtheorem{dfn}[prop]{Definition}   
\newtheorem{theo}[prop]{Theorem}
\newtheorem{conj}[prop]{Conjecture}
 
\newtheorem{coro}[prop]{Corollary} 
\newtheorem{lem}[prop]{Lemma}   
\newtheorem{exam}[prop]{Example}
\newtheorem{quest}[prop]{Question}
\newtheorem{problem}[prop]{Problem}

\title{Lectures on topology of complements and fundamental groups}
\author{A.Libgober\\
\small Department of Mathematics \\
\small University of Illinois at Chicago\\
\small 851 S.Morgan Str. Chicago, Illinois, 60607 \\
\small e-mail: libgober@math.uic.edu \\}

\begin{document}

\date{}

\maketitle

\begin{abstract}{This is an introduction to the 
topology of the complement to plane curves and hypersurfaces
in the projective space and based on the lectures
given in Lumini in February and in ICTP (Trieste) in August 
of 2005. We discuss key problems concerning the  
families of singular curves, the one variable Alexander polynomials 
and the orders of the homotopy groups of the complements to 
hypersurfaces with isolated singularities. We also discuss 
multivariable generalizations of these invariants and the Hodge
theory of infinite abelian covers used in calculations of multivariable
invariants. A historical overview is included as the opening section}
\end{abstract}


\section{Introduction}

Study of the topology of plane algebraic curves is a very old subject.
In fact, its problems come up naturally after 
the very first definitions in a basic course on algebraic curves. 
And yet, the answers obtained so far, are often elusive and incomplete.
If $C$ is an algebraic curve in a complex projective plane 
$\PP^2$, what is the fundamental 
group of $\PP^2-C$? Which properties of $C$ affect complexity of 
this group? For which groups $G$ there exist $C$ such the 
$G$ is the fundamental groups of the complement to $C$?  
When two curves are isotopic is appropriate sense, so that 
complements stay unchanged during such isotopies? 
What are the invariants of such isotopies?
These are obvious questions and a lot is known 
about them but complete or even   satisfactory answers still are out 
of reach. Below I want to describe 
some recent developments and I hope that this can serve as 
an introduction to these ideas and methods.

Perhaps the real 
beginning of this subject should be credited to Enriques, though 
some important work on construction of interesting singular curves and 
the numerology i.e. calculations of the number of singular 
points of a given type etc. started much earlier.
For example in early 19 century 
Pl\"ucker  discovered important formulas related 
the degree, the number of nodes and cusps of a curve to 
similar invariants of the dual 
curve, people from Newton to Puiseaux and beyond 
developed methods for analyzing the 
singular points of plane curves and already Newton classified types of 
singular cubics.
Lefschetz (\cite{lefsh1913}) used 
Pl\"ucker's work to obtain first non trivial information on how many nodes
and cusps a plane curve of a given degree can have 
(problem which still remains largely unresolved).

In the end of 19th century, undoubtedly influenced by Picard and Severi
works on the topology of complex surfaces,
Enriques initiated a program to extend
Riemann and Hurwitz results on multivalued functions, or in a more modern 
terminology covering spaces of Riemann surfaces, to higher dimensions
(cf. \cite{enriques}).
 According to Riemann, a  multivalued function in one variable, 
(e.g $w=\sqrt z$ or more generally a solution  
to the equation $w^n+a_1(z)w^{n-1}+...+a_n(z)$ where $a_i(z)$ are
single valued holomorphic functions of $z$)
is specified by the following data. Firstly, the collection of its 
ramification points 
$B \subset \CC \subset \PP^1$, secondly the number $d$  
of values of the 
multivalued function and finally the monodromy representation 
$\pi_1(\PP^1-B) \rightarrow \Sigma_d$ of the fundamental group into 
the symmetric group on $d$ letters.  What makes 
Riemann's approach very effective for the 
description of multivalued functions is the fact that 
the fundamental group in question is always
a free group since the ramification locus is 
just a collection of points in $\PP^1$,
and therefore the whole multivalued function is specified 
by ramification locus $B$ and 
assignment of arbitrary permutations $\sigma_1,...,\sigma_{{\rm Card}(B)}$
 in the symmetric group $\Sigma_d$ on $d$ letters
to the generators of $\pi_1(\PP^1-B)$ 
with the only restriction  
$\sigma_1 \cdot \cdot \cdot \sigma_{{\rm Card}(B)}=id$.

It was realized by Enriques (and others; a rather complete 
account of the work before mid 1930's is given by Zariski in his 
seminal book \cite{Zariskibook})  that similar description 
of multivalued functions of several variables is still 
valid but also that in higher
dimensions such a result is much less efficient since 
$\sigma_1,...,\sigma_{{\rm Card}(B)}$ must satisfy additional relations. 
For example, any algebraic 
curve in $\PP^2$ can  be a branching curve of a multivalued function but 
one cannot assign arbitrary elements of $\Sigma_d$ to generators 
of $\pi_1(\PP^2-B)$ since this group is almost never free. Rather, 
the permutations should satisfy certain compatibility conditions
(one should note that the very idea of
 the fundamental group did not completely 
crystallized around the time of the work of Enriques and
therefore his statements are much less straightforward than presented here). 
Enriques described this conditions very explicitly and this,
in modern terms, amounts to calculation of the quotient 
of the fundamental group
by the intersection of subgroups of finite index in terms 
of geometric generators
(those discussed in section  \ref{homotopygroupsviapencils}; note that 
it is still unknown that this intersection trivial i.e. the 
fundamental group is residually finite; cf. section 
\ref{fundgroupcompl}).
For example, if the branching curve has degree $d$ and is {\it non singular}, 
the fundamental group of the complement is cyclic of order $d$, but, at 
the same time, the number of geometric generators 
is $d$ (cf. section \ref{commfund}). 
Therefore, firstly, one can assign to one geometric 
generator only a permutation of order $d$ in $\Sigma_d$ and, 
secondly, the assignment to the 
rest of the generators is determined by the first choice. 
O.Zariski (after arrival to US and 
visiting Princeton where  
Lefschetz and Alexander were working at the time) realized that the 
fundamental group of the complement is the central object in this theory and 
introduced many new at that time ideas even in 
context of similar problems in the knot theory. 
He showed how subtle the question 
of the fundamental group can be: not that is does 
depend on degree of the branching 
curve, as is the case for multivalued functions in one 
variable, but even knowing the number of nodes and cusps is not sufficient.
He proved that a curve with 6 cusps can have as 
the fundamental group the cyclic group
$\ZZ_6$ or the free product $\ZZ_2 * \ZZ_3=PSL_2(\ZZ)$. 
He also showed that such sextics can be distinguished 
by a geometric condition: in the first case the 
cusps must be in general position i.e. not 
to belong to a curve of degree 2 and in the 
second case they must belong to a conic. 
Zariski also used many technical ideas just 
appearing at the time in topology e.g. 
studying the homology of cyclic covers (which in knot theory can be traced to 
Alexander and Reidemeister \cite{alexander}, \cite{reidemcovers}).
Systematic study of the branched coverings using the theory 
of adjoints (cf. section \ref{adjointsection}) allowed him to 
relate the homology of branched covers to the superabundances of linear 
systems defined by the cusps (cf \cite{irreg}). He found a close relationship 
between the fundamental groups of the complements and braid
groups by considering the dual curves for
nodal rational and elliptic curves. One of the tools 
was his celebrated theorem on fundamental groups of hyperplane sections
extending Lefschetz homological results.
In the context of branched covering Zariski 
even obtained expressions close to Alexander polynomial 
(cf. \cite{zariskicyclic})
 as was noticed by D.Mumford (cf. \cite{Zariskibook}) 
and which was the basis of his questions 
about the role of Alexander polynomial in algebraic geometry (*).
\footnote{${}^*${These questions were answered later in author's papers 
\cite{duke} \cite{arcata1} and further
extended in \cite{annals}, \cite{charvar} \cite{sapporo} (see 
references to other related works in these papers).}}

After 1937 Zariski abruptly changed the scope of his interests and 
turned to ambitious project of reconstructing algebraic geometry 
on firm foundations of commutative algebra. Some of his students,
however, continue to develop this subject, cf. \cite{turpin}, 
\cite{lehr}; much later, but in a similar spirit, M.Oka (\cite{oka76}) 
generalized Zariski calculation of the fundamental 
group of the complement to sextic with six cusps on conic by 
proving that for the curve $C$ given by the equation $(x^p+y^p)^q+
(y^q+z^q)^p=0, \ \ gcd(p,q)=1$ 
one has $\pi_1(\PP^2-C)=\ZZ_p*\ZZ_q$. The study of the topology 
did continue mostly in the works of O.Chisini and his students 
(\cite{chisini}) who initiated use of braids for the study of the 
fundamental groups and covering spaces. Abhyankar (cf. \cite{abhyn}), who 
studied with Zariski in Harvard in the 50s,
was investigating the fundamental groups, 
and in particular obtained important results of 
the fundamental groups of the complements,
but the main focus was on algebraization of the fundamental groups.

One of the driving problem in the study of the fundamental groups
in the 60s  and 70s
was the problem of commutativity for fundamental groups of the 
complements to curves with nodes only. Severi (\cite{severi}) outlined
an argument which, as later was realized, is incomplete.
It was based on an assertion that the variety of plane curves of fixed 
degree with a fixed number of nodes is irreducible. 
Zariski repeat Severi's argument in \cite{Zariskibook} but return to this 
issue much later (cf. \cite{irred}).
Severi's statement 
eventually was confirmed by J.Harris (\cite{harris}). 
A direct algebraic proof of commutativity was found 
by W.Fulton (cf \cite{fulton})
using Abhyankar's work and shortly after that a topological 
argument was given by P.Deligne.(cf. \cite{zarconjdeligne}) .
A little bit later, M.Nori (cf. \cite{nori}) 
clarified these results further
by obtaining conditions for commutativity of the fundamental group 
of the complement to curves on arbitrary surfaces, in this 
respect continuing the work of Abhyankar (cf. \cite{nori})

In the 70's the problems about fundamental groups of the 
complements were mentioned infrequently.
Mumford in already quoted appendix to \cite{Zariskibook}
also raised the problem of investigating of the quotient
$G'/G''$ for the fundamental groups of the complements.
In introduction to the volume III of collected papers by 
Zariski, containing the papers on the topology of the 
complements, Artin and Mazur, after discussing Zariksi's 
study of cyclic multiple planes, note:

\bigskip
{\it "Also, as far as the editors are aware, there has been no further
progress in the delicate study of cyclic multiple planes for 
general $d$. There are many tantalizing questions here-there
are even a number of less delicate topological issues to sort out.
For example, for irreducible plane curve $C$ with arbitrary singularities
can one give some reasonable sufficient conditions for regularity of 
$H_d$ in terms of zeros of ''local Alexander polynomials''- that is, 
the Alexander polynomials of the knots associated with singularities of 
$C$?''}
\bigskip

The answers to these questions were obtained in 
author's papers \cite{duke} and \cite{arcata1}
which depend on Milnor's work on the Alexander polynomials and the 
infinite cyclic covers (cf. \cite{infmilnor}). If $G=\pi_1(\CC^2-C)$ 
one has $G/G'=H_1(\CC^2-C)=\ZZ^r$ where $r$ is the number of irreducible 
components (cf. \ref{homologysection}). If $C$ is irreducible one has the 
exact sequence:

$$0 \rightarrow G'/G'' \rightarrow G/G'' \rightarrow \ZZ \rightarrow 0$$

This sequence 
defines the action of $\ZZ$ on 
$G'/G''$ which after Hurewicz identification 
of $G'/G''$ with the homology of infinite cyclic cover coincides 
with the action induced by the action of  
 the group $\ZZ$ of covering transformations on the universal cyclic cover.
The advantage of replacing projective curve by affine is that 
in affine case one has infinite tower of 
covering spaces while in projective 
case the degree of the cover must divide the degree of the curve. 
On the other hand, if the line at infinity is transversal to 
a projective curve, the group of affine curve is just a central extension 
of the projective one (in non-transversal case the relation is much 
more subtle). It is shown in \cite{duke} that $G'/G'' \otimes \QQ$, 
as a module over the group ring of $\ZZ$ i.e. the ring $\QQ[t,t^{-1}]$, 
is a torsion module and hence the order of the latter $\Delta_C(t)$
is well defined (up to a unit of $\QQ[t,t^{-1}]$). This is a {\it global} 
invariant of the curve in $\CC^2$. On the other hand, with each singular 
point of $C$ is associated the link i.e. the intersection of $C$ with 
the boundary of a small ball about this singular point. As result 
one obtains a set of local Alexander polynomials $\Delta_P$ 
of all singularities $P$ of the curve $C$ (as was suggested by 
Artin and Masur). One need, however, another important 
ingredient: in \cite{duke} was introduced
the Alexander polynomial at infinity $\Delta_{\infty}$
which is the Alexander polynomial
of the link which is the intersection of $C$ with the ball 
in $\CC^2$ of a sufficiently large radius. 
The answer to the question of Artin and Mazur in the above quote 
is given by the following divisibility theorem from 
\cite{duke} for these Alexander polynomials associated with the curve: 
$$ \Delta(C) \ \vert \ \Pi_{P \in {\rm Sing}(C)} \Delta_P(C)$$
\begin{equation}\label{introdivis}
 \Delta(C) \ \vert \ \Delta_{\infty}(C)
\end{equation}
and the theorem which expresses the homology of cyclic covers
in terms of Alexander polynomials (cf. \cite{duke} 
and theorem \ref{covershomology} below).
For example, for sextic curves with 6 cusps, which Zariski was considering 
in \cite{zariski29}, the Alexander polynomial $\Delta_C(t)$ 
is equal to $t^2-t+1$ or $1$
depending on wether six cusps are or are not on conic and both divide 
$(t^2-t+1)^6$ and $(t^6-1)^4(t-1)$ (which are the product of 
local Alexander polynomials and the Alexander polynomial at infinity 
respectively). 
The divisibility relation \cite{duke} 
also contains restrictions on the fundamental groups: $G'/G'' \otimes \QQ$
is trivial if $\Delta_{\infty}$ and $\Pi_{P \in {\rm Sing}(C)} \Delta_P(C)$
are relatively prime. For example if cusps are the only singularities, 
then $G'/G'' \otimes \QQ=0$ unless the degree is divisible by $6$.
The regularity condition, which were mentioned by Artin-Masur is 
the following: if none of the roots of local Alexander polynomial 
is a root of unity of degree $d$ the degree $n$ then the cyclic cover 
$H_n$ or degree $n$ is regular (cf. \cite{duke}).

The work \cite{duke} is topological and many of the results were 
extended to differential 
category (cf. \cite{lomonaco}). Dependence of Alexander polynomial
on position of singularities, 
containing generalization of several 
Zariski's calculation, was shown  in \cite{arcata1}.
As in \cite{irreg}, the irregularity of cyclic multiple plane is obtained
in terms of superabundances of certain linear systems associated with 
the cusps but for singularities more complicated than cusps the systems
are specified by more subtle geometric conditions: the local equations
of the linear systems responsible for irregularity of 
cyclic branched 
covers must belong to certain ideals called 
in the {\it ideals of quasiadjunction}. Later, these ideals appeared 
many other contexts and often are called {\it multiplier ideals} 
(cf. \cite{lazarsfeld}). Other important numerical invariants 
of plane singularities introduced in \cite{arcata1} were identified 
in \cite{loeservac} with the part of the spectrum introduced in 70's 
by Arnold and Steenbrink (cf \cite{arcata}). The work \cite{esnault} 
also related 
the irregularity of multiple planes to the position of singularities 
and these ideas rely on vanishing 
theorems which
later lead to much better understanding of those (cf. \cite{esnaultvan}): 
a key development in algebraic geometry in 90s.

In early 80s, about the time when described above work on Alexander 
polynomials appeared, there was another important development 
in the study of plane singular curves. B.Moishezon 
initiated program for describing the topology of algebraic
surfaces in terms of branching curve in $\PP^2$. This is a special 
class of curves and these curves belong to the class of curves 
having nodes and cusps as 
 singularities. If one starts with a projective surface, considers 
a pluricanonial embedding using a fixed multiplicity of the canonical 
class, and then uses a generic projection the branching  
curve in $\PP^2$ becomes an invariant of the deformation type of the surface
(the fact that one does not need the monodromy representation 
into symmetric group was conjectured by Chisini and subsequently 
proven in (\cite{kul}). Moishezon's first calculations deal 
with the branching curves of generic projections of non singular
surfaces in $\PP^3$. If the degree $d$ of a surface is 3, one obtains 
as the branching curve Zariski's 
sextic given by the equation $f_2^3+f_3^2=0$. For surfaces of arbitrary $d$ 
Moishezon obtains, as the fundamental groups of the 
complements to the branching curves, 
the quotients of the Artin's braid groups by the centers
(which for $d=3$ gives $PSL_2(\ZZ)$). Moishezon's important idea was that the
primary invariant is not the fundamental group but rather 
the braid monodromy which implicitly is present in van Kampen's method 
of calculation of the fundamental group (Moishezon was unaware 
of Chisini's work \cite{chisini} until he completed \cite{moishezon}).
In this vein, the author showed that the braid monodromy
defines not just the fundamental group but also the homotopy 
type (cf. \cite{homotopytype}, and further works by M.Teicher
 cf. \cite{teicher} )
Later Moishezon continued this work jointly with M.Teicher.
Methods of braid monodromy recently found applications in 
symplectic geometry (cf. \cite{symplectic}). More recently
Teicher and her students continued systematics study of the 
braid monodromy and the fundamental groups of the complements 
to the branching curves of generic projections and arrangements of 
lines.

In the late 80's the work started on a generalization of the 
the theory of complements to singular curves to higher dimensions.
The case of hypersurfaces with isolated singularities it turns out 
remarkably similar to the case of curves. In \cite{annals}, 
the author showed that for $n>1$ the role of Alexander 
polynomial plays the order of the homotopy group $\pi_n(\CC^{n+1}-V)$
considered as the module over $\pi_1(\CC^{n+1}-V)=\ZZ$. The point 
is that this homotopy 
group can be canonically identified with the homology 
$H_n(\widetilde {\CC^{n+1}-V},\ZZ)$ of the infinite cyclic cover 
of the complement. The divisibility relations (\ref{introdivis}) 
extends to the order of the homotopy groups and examples of hypersurfaces
with non trivial homotopy appears as a natural generalizations
of Zariski's sextics. For example $\pi_2(\CC^3-V)=\ne 0$ 
for $V$ given in $\PP^3$ by the equation: $f^2_{21}+f^3_{14}+f^7_5=0$
($f_n$ generic form of degree $n$ in four variables).
Analytic theory developed by the author in 
\cite{arcata1} also was extended to higher dimension
in \cite{homotopy}
by introducing the mixed Hodge structure on the homotopy group
and by relating one of the Hodge components of the homotopy group to 
the superabundance of the linear systems defined by the singularities of 
hypersurface.

In the 90's were obtained the first results on a multivariable generalization 
of the Alexander invariants (cf. \cite{topandappl}). 
The theory of the multivariable Alexander 
polynomials of links, due to R.Fox, depends on a very special 
feature of the link groups: the first Fitting ideal of the 
Alexander module is ``almost'' principal. The fundamental groups of 
the complement to reducible algebraic 
curves in $\CC^2$ are similar to the link groups in the sense that 
both have surjections on $\ZZ^r (r >1)$. However for 
algebraic curves the first Fitting ideal of the Alexander module is 
far from principal and as result one cannot define a
multivariable Alexander polynomial in a meaningful way. The puzzle
of existence of multivariable invariants of algebraic curves 
got resolved by introducing set of zeros of the polynomials in
the Fitting ideals of the Alexander modules in author's 
paper \cite{topandappl}. In the case of one 
variable Alexander polynomials no information get lost by  
replacement the Alexander polynomial by its set of zeros
(at least for curves in $\PP^2$ for which the Alexander module 
is semisimple)
but for reducible curves zero sets provide a non trivial 
and very interesting invariant. 

Applications followed shortly.
In \cite{eko} the characteristic varieties of a 
group were related to the cohomology of local systems which followed  
the study of polynomial periodicity of Betti numbers 
of branched covering spaces (cf. \cite{hironakagrowth}). 
For the curves for which all 
components have degree 1 i.e. arrangements of lines the 
components of characteristic varieties were related to the 
cohomology algebra of the complement (cf. \cite{cohensuciu}).
The calculation of the homology of abelian covers 
constructed by Hirzebruch and which have universal covers 
biholomorphic to the ball did fall naturally in the general scheme
valid for arbitrary arrangements and covers (cf. \cite{charvar}).
Analytic (rather than topological) theory 
was developed in \cite{charvar} and characteristic 
varieties were expressed in terms of superabundances of the linear systems.
Essential in this calculation were the results in \cite{arapura}, 
on the structure of the jumping loci for the cohomology of local systems.
They represent an extension to quasiprojective varieties of the results
of Green-Lazarsfeld, Beauville, Catanese, Simpson, Deligne 
 and others and which 
asserts that the jumping loci for the cohomology of local systems
are cosets of certain subgroups of the group of characters 
of the fundamental group.

During late 90's, the study of the topology of plane algebraic curves 
became much more active area of research. Many new examples 
of Zariski pairs due to E.Artal-Bartolo 
and collaborators and M.Oka showed how common is the phenomenon 
of curves having different equisingular isotopy type 
with the same local data. Many new calculations were carried out 
of the fundamental groups of the complements by M.Teicher's school 
which finally led to a general conjecture on the structure 
of the fundamental groups of the branching curves of generic 
projections (cf. \cite{teicher}). 
Interactions with combinatorics of arrangements 
were important and lead to at least conjectural description 
of the characteristic varieties and much stronger vanishing 
for the cohomology of local systems than were available 
earlier (cf.\cite{sergey},\cite{charvar},\cite{cohensuciu}).
Connections with symplectic topology should be noted (cf. 
\cite{symplectic}). 
There was further progress in the study of the complements
in higher dimensions on generalizations of Zariski-van Kampen's
theorem (cf. \cite{annals},\cite{denisme},\cite{eyral},\cite{tibar}).
Nevertheless, despite 
tremendous progress, since the first works by Enriques, Zariski
and van Kampen, many problems still remains open and complete 
understanding of the topology of the complements to curves
and hypersurfaces still out of reach. 

In the text below 
we outlined some of the problems which resolution may clarify 
substantially the situation. Exposition is very elementary in the beginning
describing motivation for the study in the following sections.
In the later parts a reader will need more and more rely on
material covered in standard courses in algebraic geometry. Moreover, 
some familiarity with the mixed Hodge theory is needed in the 
last sections. Textbook \cite{dimca} also is a good reference
to the background and other material discussed below.
Most of the material appeared already in the literature some 
time ago but some 
results appear to be new.

I want to thank J.P. Brasselet, D.Cheniot, J.Damon, M.Oka, A.Pichon
D.Trotman,  N.Dutertre and C. Murolo
who organized the conferences in Lumini and Trieste and the 
opportunity to present this beautiful area of mathematics.
I am also want to thank L.Kauffman for a discussion of the
history of polynomial invariants in the knot theory.

\section{Fundamental groups of the complement}

\subsection{Problem of classification up to isotopy.}  

\bigskip

\subsubsection{Stratification of the discriminant}

Classically, many problems in the topology of 
plane curves and hypersurfaces were 
rooted in the study of families of curves and attempts of 
some kind of classification of curves and 
 hypersurfaces (cf. \cite{zariski29}). 
We shall start by discussion in what kind classification and in which 
sense one may expect.
\par Hypersurfaces of a fixed degree $d$ in $\PP^n$ are parameterized by 
$\PP^{{{n+d} \choose d}-1}$ by assigning to a defining equation 
the collection of coefficients of its monomials (in some fixed order). 
The discriminant $Disc(n,d)$ is the hypersurface in $\PP^{{{n+d} \choose d}-1}$
consisting of the points corresponding to singular hypersurfaces.
It has singularities in codimension one.
An interesting problem is to understand the stratification 
of the discriminant hypersurfaces $Disc(n,d)$. By this we mean to 
describe the singular 
locus of the discriminant hypersurface having codimension one in
$\PP^{{{n+d} \choose d}-1}$, then the singular locus of 
singular locus (having codimension 2) and so on. 
More precisely, we consider the universal hypersurface of degree $d$ i.e. 
${\cal V} \subset \PP^n \times \PP^{{{n+d} \choose d}-1}$ consisting 
of pairs $(P,V)$ such that $P \in V$. $Disc(n,d)$ is the image 
of the critical set of the projection on the second factor 
and its preimage in $\cal V$ is the universal singular hypersurface.
The critical set of the projection on the second factor is 
the singular set $Sing(n,d)$ of the universal singular hypersurface.
Then we consider the critical set of the restriction of the 
projection on the non singular part of $Sing(n,d)$. 
On a codimension one subset, the rank
of projection on $Disc(n,d)$ drops and so on.
With such a definition, Thom's isotopy theorem yields that 
the hypersurfaces belonging to each stratum are equisingular so 
the strata represent equisingular families of hypersurfaces. Note that 
the subset in $\PP^{{{n+d} \choose d}-1}$ 
parameterizing equisingular hypersurfaces is singular in general
(cf.\cite{wahl}).

The case $n=1$ is 
already very interesting and non trivial.
The discriminant consist of homogeneous polynomials $\prod_i(\alpha_i u
-\beta_iv)$ in two variables $u,v$ having multiple roots,
i.e factors such that $(\alpha_i,\beta_i)$ and $(\alpha_j,\beta_j)$
satisfy $det\vert  \matrix{\alpha_i & \beta_i \cr \alpha_i & \beta_i}
\vert=0 $. The strata correspond to partitions of $d$, i.e. the 
conjugacy classes of the symmetric group $\Sigma_d$. 
A lot is known about the geometry of these strata, for example the 
degrees of their closures as well as other algebro-geometric information.
Cases with $n>1$ are much more complicated. 
Many pieces of information are known.  For example, in the 
case $n=2$ the degrees of the stratum corresponding to 
rational nodal curves have the interpretation as
Gromov-Witten invariants of a projective plane and as such 
satisfy beautiful recurrence relations (cf. \cite{KM}). 
Indeed, the dimension 
of this stratum is $3d-1$ where $d$ is the degree of the curves
(i.e. ${{(d+1)(d+2)} \over 2}-{{(d-1)(d-2)} \over 2}$) so the degree
of the corresponding stratum 
is the number of nodal curves of degree $d$ passing through generic
$3d-1$-points.  The degrees of strata of nodal curves are subject
to a conjecture of G\"ottsche discussed, for example, in \cite{kleiman}.

\subsubsection{Classification of quadrics and cubics. Local type.}
Another class of discriminants which is 
well understood is  case $d=2$. Each stratum 
correspond to the quadrics of fixed rank. In particular 
each stratum is a determinantal variety. 

Classification of plane cubics goes back to Newton. Codimension 
one stratum consists of cubic curves with one node. It has the 
degree equal to 12. There are two strata having codimension 2. 
One consists of curves with one cusp and another formed by 
the reducible curves having as the components a nonsingular quadric 
and a non tangent to it line.
The rest of the strata correspond to reducible curves and 
each is determined by strata of curves of lower degree and 
the mutual position.
The strata of codimension three are: unions 
of a non singular quadric and a tangent line (in the 
closure of both strata of codimension 2) and the union of 
three lines in general position. Note that each of these strata 
is described by the local type of singularities: 
the number of nodes, cusps, tacnodes etc. A definition of the 
local type is the following:

\begin{dfn}\label{localtype} 
Two reduced curves $C$ and $C'$ (of the same degree) 
have the same local type if there
is a one to one correspondence between their singular points such the each 
pair of corresponding singular points $P$ and $P'$ have neighborhoods
$B_{\epsilon}$ and $B'_{\epsilon'}$ and homeomorphisms $\phi_P: B_{\epsilon}
\rightarrow B'_{\epsilon'}$ such that $\phi_P(C \cap B_{\epsilon})=
C' \cap B'_{\epsilon'}$. 

\noindent Two possibly non reduced curves have the same local type if: 

\noindent (a) corresponding reduced curves have the same local type 

\noindent (b) there are one to one correspondences between the components and 
singular points such that corresponding components have the same degrees and 
multiplicities and corresponding singular points belong to the corresponding components.  
\end{dfn}

\subsubsection{Examples with disconnected strata.} The classification 
of strata of curves of degree 4 provides the first example 
when the local type of singularities (in the sense of the first 
part of the definition \ref{localtype})
yields the strata with several connected components. The quartics  
with three nodes have two types: firstly the irreducible ones 
and quartics which are the unions of a non singular cubic and 
a generic line. The strata are distinguished by a global property.

For each degree there are finitely many irreducible families 
of plane curves having the same local type.
\begin{problem} Find discrete invariants of families of 
curves having the same local type.
\end{problem}
This problem is similar to the problem of classification 
of knots in $S^3$. Thom's isotopy theorem implies that the 
curves (or hypersurfaces) in a {\it connected} equisingular family 
are isotopic and hence have 
diffeomorphic complements. 
The main tool in the study of knots is the 
fundamental group of the complement which is one of the 
reasons suggesting to look at $\pi_1(\PP^2-C)$ or also into 
$\pi_k(\PP^n-V)$ with $k>1$ in the case of hypersurfaces of higher 
dimensions. 

\subsection{Fundamental groups of the complements}\label{fundgroupcompl}

Though the classification problem of the strata of 
the discriminant leads to the 
fundamental groups of the complements as potentially an important 
invariant there are many other 
reasons for looking at the fundamental groups.
One is that the fundamental groups of the complements to 
hypersurfaces controls the covers of projective space and any projective
algebraic variety having the dimension $n$ 
is a branched covering space of $\PP^n$.

Linear representations of the fundamental groups appear as the
monodromy representations of differential equations 
and this correspondence is a subject of the Riemann-Hilbert problem.
For example monodromy representation of KZ equation yields 
an interesting representation of the pure braid group closely 
related to the discriminant $D(1,d)$.

Each of these ``applications'' lead to concrete problems about the 
fundamental groups. For example, the use for a study of the coverings
suggests the following. In the above presentation 
of algebraic varieties as the covering spaces the degree 
of a cover is always finite. So the coverings are determined 
already by the quotient by the intersection of all subgroups of 
finite index. The problem is if this intersection is the identity group
or in other words if 
the fundamental group of the complement to an algebraic hypersurface
is {\it residually finite}.  Alternatively this the the question on 
whether the map $\pi_1 \rightarrow 
\pi_1^{alg}$ into the algebraic fundamental group is injective.
Note that the fundamental group of an algebraic variety does 
not have to be residually finite (D.Toledo).
In general the problem of finding the properties of the fundamental groups 
of the complements or characterizing the algebraic structure of these
groups is one of the central and the most difficult 
problems in algebraic geometry.

\subsubsection{Homology of the complements.}\label{homologysection}
 An easily available 
information about the fundamental groups $\pi_1(\PP^{n+1}-V)$ 
comes from calculation of the  
homology $H_1(\PP^{n+1}-V)$ which, by Hurewicz theorem, is the 
quotient of the fundamental group by its commutator. 
Here is an answer:

\begin{prop}\label{homology}
 Let $V$ is the union of irreducible components $V_1,..,V_r$ 
having the degrees $d_1,...,d_r$. Then 
$H_1(\PP^{n+1}-V,\ZZ)=\ZZ^{r}/(d_1,..,d_r)$.
\end{prop}

For example, if $g.c.d.(d_1,...,d_r)=1$ then the homology group 
is torsion free. This is the case when one of components has the degree 
equal to 1 or in other words for the complements to affine hypersurfaces.

\subsubsection{Examples of calculations of the fundamental groups.}

\bigskip

In the last twenty years quite a few calculations of the 
fundamental groups were made. For example, as mentioned in 
introduction,  Moishezon-Teicher calculated 
the fundamental groups to the branching curves of generic projections
of many algebraic surfaces 
(cf. \cite{teicher}). Oka calculated the fundamental 
groups of the complements 
to many curves having low degree, in particular to various 
classes of curves of degree 6 (cf. \cite{oka}). 
Many calculations were carried out by Artal-Carmona-Cogolludo 
(cf. \cite{artal}).
These techniques I will discuss in a later chapter,  but here I want 
to explain some short and elegant calculations made by  
Zariski 80 years back.

\begin{prop}\label{ratnodal} 
Let $\hat C_d$ be a curve dual to the rational nodal curve $C_d$
having the degree $d$ (the degree of $\hat C_d$ is equal to $2(d-1)$, 
it has $3d-6$ cusps  and $2(d-2)(d-3)$ nodes). The group 
$\pi_1(\PP^2-\hat C_d)$ 
is isomorphic to the braid group of sphere on $d$ strings. 
In particular the fundamental group of the complement to the 
quartic with 3 cusps is a non abelian group having order 12.
\end{prop}

Indeed, $C_d$ is a generic projection on $\PP^2$ of a rational normal curve
$C$ in $\PP^d$ and dual to $C_d$ curve is a plane 
section of the hypersurface $\hat C$ in $\PP^d$ dual to $C$ by 
a plane $H$. The complement to 
this hypersurface consists of hyperplanes in $\PP^d$ 
intersecting $C$ in $d$ distinct points which can be chosen arbitrary. 
Hence the space of based loops in this complement is 
identified with the braids of $\PP^1(\CC)=S^2$. Finally 
the isomorphism $\pi_1(\PP^2-C_d)=B_d(S^2)$ follows from 
Lefschetz hyperplane section theorem applied to embedding 
of the complement in $H$ into the complement in $\PP^d$.
In the case $d=3$, the pure braid group of sphere can be identified
with $\pi_1(PGL_1(\CC))=\ZZ_2$ and hence one has the exact sequence:
$1 \rightarrow \ZZ_2 \rightarrow B(S^2) \rightarrow S_3 \rightarrow 1$.

\subsubsection{Alexander invariants of the fundamental groups.}
\label{alexsinvariants}

\bigskip

Since the problem of characterization and understanding the
fundamental group is very complicated it is reasonable to 
try to rather understand some invariants of the  fundamental groups.
An accessible and interesting invariant is the Alexander invariant of 
a group.

Let $G$ be arbitrary group together with a surjective homomorphism 
$\phi: G \rightarrow \ZZ^r$. Let ${\rm Ker} \phi =K$ and let $K'=[K,K]$
be the commutator. If $\phi$ is the abelianization $G \rightarrow G/G'$ 
then $K=G', K'=G''=[G',G']$. 
We have: 

\begin{equation}
0 \rightarrow K/K' \rightarrow G/K' \rightarrow \ZZ^r \rightarrow 0
\end{equation}

In particular $K/K'$ receives the action of $\ZZ^r$ and 
hence $K/K'$ becomes the module over the group ring of $\ZZ^r$.
This module is called the Alexander 
invariant of the pair $(G,\phi)$. In the 
case when $\phi$ is the abelianization one obtains 
an invariant depending on the group $G$ only. 
It is denoted below as $A(G,\phi)$ or
if $\phi$ is the abelianization as $A(G)$.

This definition can be interpreted geometrically. If $X_G$ is a 
CW-complex having $G$ as its fundamental group then  
the homomorphism $\phi$ defines a covering space $\tilde X_{G,\phi}$.
One has $\pi_1(\tilde X_{G,\phi})=K$ and $K/K'=H_1(\tilde X_{G,\phi},\ZZ)$.
The action of $\ZZ$ corresponds to the action of the group 
$\ZZ$ of deck transformations on $\tilde X_{G,\phi}$.

For perfect groups, i.e. such that 
$G=G'$, this invariant is 
trivial (since $r=0$ is the only possibility), but 
 since the fundamental groups of the complements in $\PP^{n+1}$
are perfect only if the hypersurface is the hyperplane (cf. \ref{homology})
for them 
the Alexander invariant is always interesting.

There is an algorithmic procedure for calculation of 
the Alexander modules due to R.Fox (''Fox calculus'') (cf. \cite{fox}).

Let $G$ be a finitely generated, finitely presented group
 i.e. one has a surjective map $\Phi: F_s \rightarrow G$ 
of the free group $F_s$ on $s$ 
generators, $x_1,...x_s$ with the kernel being the normal closure
of a finite set of elements given by the words $R_1,... R_N$ in
$F_s$. Consider the 
maps of the group rings: ${\partial \over {\partial x_j}}:
 \ZZ[F_s] \rightarrow \ZZ[F_s]$ uniquely 
specified by the conditions: 
\begin{equation}
{{\partial (uv)} \over {\partial x_j}}={{\partial u} \over {\partial x_j}}a(v)+
u{{\partial u} \over {\partial x_j}};  \ \ \ \ \ \ \ \ \ \ \
{{\partial x_i} \over {\partial x_j}}=\delta_{i,j}
\end{equation} 
where $a: \ZZ[F_s] \rightarrow \ZZ$ is the augmentation surjection.
Using operators $\partial \over {\partial x_j}$ 
one can define the map of  
free $\ZZ[\ZZ_r]$-modules given by the Jacobi matrix:
\begin{equation}\label{foxjacobi}
(\phi_* \circ \Phi_*{{\partial R_i} \over {\partial x_j}}):  
\ZZ[\ZZ_r]^N \rightarrow
  \ZZ[\ZZ_r]^r
\end{equation}
which entries are obtained by applying homomorphisms 
$\Phi_*: \ZZ[F_s] \rightarrow \ZZ[G]$ 
and $\phi_*: \ZZ[G] \rightarrow \ZZ[\ZZ^r]$
of group rings induced by $\Phi$ and $\phi$ respectively.
The geometric meaning of this map is the following.
With a presentation $\Phi$ one can associate the 2-complex $X_G$ 
with single $0$-cell, $r$ $1$-cells forming wedge $S^1 \vee ... \vee S^1$
of circles corresponding to the generators 
of $G$ and $N$ 2-cells attached so that the boundary of each 
is represented by the word $R_i \ \ (i=1,...,N)$ 
in $S^1 \vee ... \vee S^1$.
The covering space  $\tilde X_{G,\phi}$ corresponding to the
homomorphism of the fundamental group  has a canonical 
cell structure given by the preimages of cells in 
the above cell decomposition of $X$: each cell in $X_G$ is replaced
by cells of the same dimension 
corresponding to the elements of the covering group.
Hence we obtain the isomorphisms $C_2(\tilde X_{G,\phi})=
\ZZ[\ZZ^r]^N, C_1(\tilde X_{G,\phi})=
\ZZ[\ZZ^r]^s $. Moreover, after this identification, 
the boundary operator $\partial_2: C_2(\tilde X_{G,\phi})
 \rightarrow C_1(\tilde X_{G,\phi})$ 
becomes identified with 
the operator given by (\ref{foxjacobi}). 
Since $H_0(\tilde X_{G,\phi},\ZZ)=\ZZ$
and $C_0(\tilde X_{G,\phi},\ZZ)=\ZZ[\ZZ]$ we have the isomorphism 
${\rm Im} \partial_1={\rm Ker}C_0(\tilde X_{G,\phi},\ZZ) \rightarrow \ZZ=
I_{\ZZ[\ZZ^r]}$ where $I_{\ZZ[\ZZ^r]}$ is the augmentation ideal of 
the group ring. Hence, (\ref{foxjacobi}) determines the presentation 
of the module very closely related to $H_1(\tilde X_{G,\phi})$.
More precisely, if $\hat H(X_{G,\phi})$ is the module having presentation 
(\ref{foxjacobi}) then  we have:

\begin{equation}
0 \rightarrow H_1(X_{G,\phi}) \rightarrow \hat H(X_{G,\phi})
\rightarrow I_{\ZZ[\ZZ^r]} \rightarrow 0
\end{equation}

For example for (affine portions of ) the curves 
in proposition \ref{ratnodal}, 
for the Alexander module $A(G)$ coinciding with $H_1(X_{G,\phi})$ 
one has (in these examples $\phi: \pi_1 \rightarrow H_1=\ZZ$ is the 
canonical homomorphism): 

\begin{equation}\label{alexbranching}
A(\pi_1(\PP^2-C_d))=\ZZ[t,t^{-1}]/(t^2-t+1) \ \ (d=4), A(\pi_1(\PP^2-C_d))=0
\ \ (d \ge 5)
\end{equation}

For the links of algebraic singularities, which all belong to the 
class of iterated torus link, the Alexander polynomial, {\it i.e. the 
order of $A(G)\otimes \QQ$ as a $\QQ[t,t^{-1}]$-module}, 
can be found 
using the data of iterations and the values of Alexander polynomial 
of for the torus knot: for the link of singularity $x^p=y^q \ \ g.c.d.(p,q)=1$
one has the following.
\begin{equation}\label{localalexformula}
\Delta(t)={{(t^{pq}-1)(t-1)} \over {(t^p-1)(t^q-1)}}
\end{equation}

Another way to calculate the Alexander polynomial is to use A'Campo 
formula for the zeta-function of the monodromy in terms of 
a resolution of the singularity (cf. \cite{EN}):
\begin{equation}
\zeta(t)=\Pi (1-t^{m_i})^{\chi(E_i^{\circ})}
\end{equation}
where $E_i$ are the exceptional curves of the resolution, $E_i^{\circ}$
is set of points in $E_i$ which are non-singular points of the exceptional 
divisor, $m_i$ is the order along $E_i$ of the pullback of the equation of the 
singularity and $\chi$ denotes the Euler characteristic. The $\zeta(t)$ 
determines the Alexander polynomial of a curve singularity via:
$\zeta(t)={{(t-1)} \over {\Delta(t)}}$. 

\subsubsection{Alexander polynomials of plane algebraic curves: 
divisibility theorems.}

There are two types of general results concerning the 
Alexander invariants of the fundamental groups $\pi_1(\CC^2-C)$.
Firstly, the Alexander polynomials of plane algebraic curves are restricted 
by the degree of the curve, by the local type of singularities and 
by position of the curve relative to the line at infinity. 
These restrictions sometimes yield triviality of the 
Alexander polynomial.
On the other hand, the Alexander polynomial is completely 
determined by the 
local type and superabundances of certain linear systems given
by the data depending on the singularities.

We shall start, with discussion of the first group of results.
Let $C$ be projective curve and $L$ be the line at infinity.
One has the linking number homomorphism: 
${\rm lk}: \CC^2-C \rightarrow \ZZ$
associating to a loop $\gamma$ in $\CC^2-C$ the (oriented) number of 
intersection points of $C$ and an immersed disk with boundary $\gamma$.
In the case when $C$ is irreducible this homomorphism 
$H_1(\CC^2-C) \rightarrow \ZZ$ already was used above. This 
defines the Alexander module and the Alexander polynomial $\Delta_C(t)$
and we shall omit mentioning the linking homomorphism used in
its definition.

To each singular point $P \in C \in \PP^2$ we associate the 
local Alexander polynomial which is the Alexander polynomial of the 
link defined as follows. In the case when $P \notin L$,the link
 is the intersection of $C$ with a sufficiently small ball
about $P$ (so that the link type is independent of the radius).
If this link has several components (i.e. $P$ has several 
branches) the Alexander polynomial again is calculated relating to 
the total linking number homomorphism. In the case when
 $P \in L$, i.e. the curve has singularities at infinity, the local 
Alexander polynomial is defined as above but $P$ considered 
as the singular 
point of $P \in C \cup L$.
Note that, as follows from definitions, the local Alexander polynomials 
can be calculated as the characteristic polynomials of the monodromy 
operators (cf. \cite{EN}, \cite{singularpoints}
 for examples and algorithms).

On the other hand, one can define the Alexander polynomial at infinity 
$\Delta_{\infty, C}$
as the Alexander polynomial of the link which is the intersection of $C$ 
with the boundary of a sufficiently small tubular neighborhood of $L$ 
in $\PP^2$ (or alternatively the sphere of a sufficiently large radius in 
$\CC^2=\PP^2-L$). For example, is $C$ is a union of $d$ lines passing through 
a point in $\PP^2$ outside of $L$ then the link at infinity is the Hopf
link with $d$ components and hence its Alexander polynomial is:

\begin{equation}\label{infinitylink}
             \Delta_{\infty,C}=(t^d-1)^{d-2}(t-1)
\end{equation} 

The same equality holds for a curve which is transversal to 
the line at infinity since there is a deformation of such a curve 
to a union of $d$ lines as above, such that transversality holds
for all curves appearing during the deformation.

With these definitions we have the following:

\begin{theo}\label{divisibility}
$$\Delta_C(t)\ \  \vert \ \ \Pi_{P \in {\rm Sing}C} \Delta_P(t)$$
$$\Delta_C(t) \ \ \vert  \ \ \Delta_{\infty,C}(t)$$

\end{theo}

Consider, for example and irreducible curve in $\PP^2$ having 
ordinary cusps (i.e. having $x^2=y^3$ as the local equation)
and nodes (local equation: $x^2=y^2$) as the only singularities. 
Then, as follows from (\ref{localalexformula}), the local 
Alexander polynomials for each singularity 
is $t^2-t+1$ (cusp) or $t-1$ (node). Moreover, it is not hard to 
show that the multiplicity of the factor $(t-1)$ is $r-1$ where 
$r$ is the number of irreducible component of $C$ (cf. \cite{duke}).
Hence we obtain:

\begin{coro} Let $C$ be an irreducible curve in $\PP^2$ having cusps
and nodes as the only singularities. Then:
  $$\Delta_C(t)=(t^2-t+1)^s$$
for some integer $s \ge 0$.
\end{coro}
Combining this corollary, the divisibility and the formula 
(\ref{infinitylink}) we obtain:

\begin{coro}Let $C$ be an irreducible curve in $\PP^2$ having cusps
and nodes as the only singularities. Then $\Delta_C(t)=1$
unless $d$ is divisible by $6$.
\end{coro}

We leave as an exercise for a reader to work out that 
$pq \not \vert d$ is a sufficient condition 
for triviality of the Alexander polynomial for an irreducible curve 
of degree $D$ with singularities locally given by $x^p=y^q$.
 
Since the curves discussed in Proposition \ref{ratnodal} 
(and also the branching curves of generic projections of non-singular 
surfaces in $\PP^3$ cf. \cite{moishezon}) have the degree $d(d-1)$ 
it follows that the Alexander polynomial is trivial if $d \equiv 2 
({\rm mod} 3)$ which explains with no calculation 
one third of equation \ref{alexbranching} (at least after 
tensoring with $\QQ$). Many additional examples of calculations 
of the Alexander polynomials can be found in \cite{oka}. 
Note, finally that it is also beneficial to consider the 
Alexander polynomials over finite 
fields ${\bf F}_p$, rather than over $\QQ$
i.e. $H_1(X_{G,\phi},{\bf F}_p)$ (cf. \cite{lomonaco}).

\subsubsection{Alexander polynomials of plane algebraic curves: 
position of singularities.}\label{alexpolynomialposition}

Now we shall discuss the dependence of the Alexander polynomial of
the positions of singularities of the curve. To this end 
we shall associate the following invariants of the 
singularities of plane curves:
rational numbers, called the constants of quasiadjunction: 
$\kappa^P_1,...,\kappa^P_{n(P)}$ corresponding to 
each point  $P$  in the set ${\rm Sing}C \subset \PP^2$ 
of singular points of $C$. 
Moreover, to each $\kappa \in \QQ$, which is a constant of 
quasiadjunction of a point $P \in {\rm Sing} C$
and each $Q \in {\rm Sing} C$,
we associate the ideal ${\cal J}_{\kappa} \subset {\cal O}_Q$
in the local ring of $Q \in \PP^2$.
($P$ and $Q$ may be distinct). This data of constants of quasiadjunction 
and the ideals in the local rings of singular points determines
the global Alexander polynomial $\Delta_C(t)$ completely (cf. 
\cite{arcata1} and 
theorem (\ref{positioncurves}) below).

The idea of calculation is based on the relation between the 
Alexander polynomial and the homology of cyclic covers on one side 
and the classical method of adjoints to describe holomorphic
2-forms on hypersurfaces in projective space (cf. \cite{Zariskibook}).

The relationship between the Alexander polynomial and the homology 
of cyclic branched covering is the following:

\begin{theo}\label{covershomology}
 Let $f(x,y)=0$ be the equation of a curve $C \in \CC^2$. 
Let $\tilde V_n$ be a desingularization of a compactification 
of the surface $z^n=f(x,y)$ in $\CC^3$. If $A(\pi_1(\CC^2-C)) \otimes \QQ=
\oplus \QQ[t,t^{-1}]/(\delta_i(t))$ is the cyclic decomposition
of the Alexander module of $C$ (i.e. $\Delta_C(t) =\Pi_i \delta_i(t))$
then ${\rm rk} H_1(V_n,\QQ)$ is equal to the sum over $i$ of the 
number of common roots of $t^n-1$ and $\delta_i(t)$. 
If the line at infinity is transversal to $C$ then the Alexander module is 
semisimple and the dimension of the $\omega_n$-eigenspace of a generator 
of the Galois group $\ZZ_n$ acting on $H_1(V_n,\CC)$ ($\omega_n$ 
is a root of unity of degree $n$) is equal to the multiplicity of 
$\omega_n$ as a root of the Alexander polynomial.
\end{theo}

Note that the first Betti number of a non-singular
projective algebraic surface is a birational invariant and hence the first 
Betti number of a resolution of a compactification is a well 
defined invariant of an affine surface $z^n=f(x,y)$. Therefore it 
also an invariant of affine curve $C$. Similar result is valid 
for branched covering of $S^3$ branched over a link: the idea 
of using covering spaces to derive invariants of knots 
goes back to Alexander and Reidemister (cf \cite{alexander}, 
\cite{reidemcovers},  \cite{zariskicyclic}, \cite{Zariskibook}). 
A consequence of this theorem is that the homology of cyclic 
covers, in the case when line at infinity is transversal to $C$ 
determine the Alexander polynomial. Another consequence, is periodicity 
of the homology of cyclic covers. In the abelian case the growth 
of the homology is polynomial periodic (cf. \cite{hironakagrowth})

The calculation of the 
homology of cyclic covers using theory of adjoints was carried out 
in \cite{irreg} (the case when $C$ has cusps and nodes), \cite{lehr} 
(the case when $C$ has singularities of the form $x^k=y^k$ or $x^k=y^{k+1}$) 
and, much later, for the curves with arbitrary singularities, in 
\cite{arcata1}. The proofs for a generalization to situation including 
hypersurfaces having arbitrary dimension is given in \cite{position}.
In fact all these proofs yields the irregularity $q={\rm dim} H^1(\tilde V_n,
{\cal O}_{\tilde V_n})={\rm dim} H^0(\tilde V_n,
{\Omega^1}_{\tilde V_n})={1 \over 2} {\rm dim} H^1(\tilde V_n,\CC)$
(and in \cite{position} the Hodge number $h^{n,0}$ for
cyclic coverings of $\PP^{n+1}$).

For details of the using this method we shall refer to \cite{position}
and section \ref{adjointsection}, 
but here we shall only remark that the adjoint ideal of a germ 
$(W,P) \in (\CC^{3},P)$
of isolated singularity at $P$ consists of germs in ${\cal O}_P$
which restriction to $W$ belongs to  
$\Phi_*(\Omega^2_{\tilde W})$ where 
$\Phi: \tilde W \rightarrow W$ is a resolution 
of singularities of $W$. If $W$ is given by an equation $F=0$, 
then the 2-forms on $W-P$ are residues of 3-forms 
${\psi(x,y,z){dx \wedge dy \wedge dz} \over F(x,y,z)}$ 
on $\CC^3$ having pole of order one along $W$ i.e. 
restrictions of 2-forms: 
\begin{equation}\label{difform}
\psi(x,y,z){dx \wedge dy } \over 
{\partial F(x,y,z) \over {\partial z}}
\end{equation}
on $W-P$.

On the other hand the 2-forms on a resolution can be described 
as the 2-forms on 
$W-P$ which can be extended over the exceptional locus of $\Phi$. 
Hence a germ $\psi(x,y,z)$ is in the adjoint ideal of $W$ if the 
pull back of the form (\ref{difform}) on resolution 
$\tilde W$ extends over the exceptional set.
Such interpretation of 2-forms on resolutions allows to relate
the dimensions of space 1-forms on $\tilde V_n$ (which is isomorphic to  
$H^1(\Omega^2_{\tilde V_n})$) to $H^1$ of certain sheaf of ideals on 
$\PP^2$ which we are going to describe.

Let $\phi(x,y)$ be a germ of a holomorphic function. Let us 
consider the function $\Xi_{\phi}(n)$ which assigns to 
a $n$ the minimal $k$ such that $z^k\phi(x,y)$ belongs to the adjoint ideal
of the singularity $z^n=f(x,y)$. 

\begin{lem} There exist $\kappa_{\phi} \in \QQ$ (also depending on 
singularity $f(x,y)$)
such that
for $\Xi_{\psi}(n)=[\kappa_{\phi}n]$ ($[.. ]$ denotes the integer part) 
\end{lem} 

The adjoint ideal of a function $F(x,y,z)$,  which is generic for its Newton
polytope, can be described as follows: a monomial 
$x^{\alpha}y^{\beta}z^{\gamma}$ 
is in the adjoint ideal of $F(x,y,z)$ if and only if the point
$(\alpha+1, \beta+1, \gamma+1)$ is inside the Newton polytope of $F(x,y,z)$
(cf. \cite{merle}). Hence if $f(x,y)=x^a+y^b$ and $\phi(x,y)=x^iy^j$
then $z^kx^iy^j$ is in the adjoint ideal of $z^n=f(x,y)=x^a+y^b$ if and only if
$(i+1)bn+(j+1)an+(k+1)ab >abn$ or 
$k+1 >n(1-(i+1){1 \over a}-(j+1){1 \over b})$. 
Therefore:
\begin{equation}
  \Xi_{x^iy^j}(n)={\rm max}([n(1-(i+1){1 \over a}-(j+1){1 \over b})],0)
\end{equation}

This construction can be used to associate to a constant $\kappa \in \QQ$ 
the following ideal 
in the local ring of the singular point of germ $f(x,y)$:

\begin{dfn} Let $\kappa \in \QQ$. The corresponding ideal of quasiadjunction
is defined as follows:
$${J}_{\kappa}=\{ \phi(x,y) \vert \kappa > \kappa_{\phi} \}$$
\end{dfn}

For example if $f(x,y)=x^2+y^3$ then:
\begin{equation}\label{quasiadjcusp}
\kappa_{x^iy^j}(n)=\cases{[{n \over 6}], & (i,j)=0 \cr 0, & $i+j \ge  1$ \cr} 
\end{equation}
and hence there is only one constant of quasiadjunction $\kappa={1 \over 6}$.
Moreover, ${J}_{1 \over 6}$ is generated by monomials such that
$i+j \ge 1$ i.e. is the maximal ideal.

Loeser and Vaqui\'e showed that the constants of quasiadjuction are precisely 
the elements of Arnold-Steenbrink spectrum of singularity $f(x,y)$
which are belong to the interval $(0,1)$. In particular 
$exp(2 \pi i \kappa_{\phi})$ are the eigenvalues of the monodromy of 
$f(x,y)=0$ and hence are the roots of the Alexander polynomial of the link of 
$f(x,y)$. After introduction of multiplier ideals it was 
soon realized that the ideals of quasiadjunction are closely related to 
multiplier ideals (cf. section \ref{multidealsection}
below). J.Kollar noticed the connection between the log-canonical 
threshold and the constants of quasiadjunction (cf. \cite{kollar} and
section \ref{multidealsection}).

Using the ideals ${J}_{\kappa}$ in the local rings of 
points in $\PP^2$, which are the singular points of a curve $C \in \PP^2$,
one defines the ideal sheaf 
\begin{equation}{\cal J}_{\kappa}={\rm Ker} {\cal O}_{\PP^2} \rightarrow
\oplus_{P \in {\rm Sing} C} {\cal O}_P/J_{\kappa,P}
\end{equation}
where $J_{\kappa,P}$ is the ideal corresponding to the singularity of $C$ 
at $P$ and the constant $\kappa$.
Using this we can calculate the Alexander polynomial as follows:

\begin{theo}\label{positioncurves} Let $C$ be a curve in $\PP^2$ 
having degree $d$ and 
let $\kappa_1, ...., \kappa_N$ be the collection of all constants 
of quasiadjunction of all singular points of $C$. Then the Alexander 
polynomial $\Delta_C(t)$ is given by:
$$\Pi_{i,d\kappa_i \in \ZZ} [(t-exp(2 \pi \sqrt{-1} \kappa_i)
(t-exp(2 \pi \sqrt{-1} \kappa_i)]^{{\rm dim} H^1(\PP^2, {\cal J}_{\kappa_i}
(d-3-d\kappa_i))}$$
\end{theo}

Note that the exponent can be written as follows:
\begin{equation}\label{defsuperab}
{{\rm dim} H^1(\PP^2, {\cal J}_{\kappa_i}
(d-3-d\kappa_i))}={{\rm dim} H^0(\PP^2, {\cal J}_{\kappa_i}
(d-3-d\kappa_i))}-\chi({\cal J}_{\kappa_i})
\end{equation}
(since $H^2(\PP^2,{\cal J}_{\kappa_i}(d-3-d\kappa)=0$).
In other words the exponent is the difference between the actual
and ``expected'' dimensions of the linear system of curves
of degree $d-3-d\kappa_i$ which local equations belong to the
ideals of quasiadjunction corresponding to the constant $\kappa$.
Therefore, (\ref{defsuperab}) is
 what is classically called the superabundance of this linear system.

As an example, let us consider the sextics with six ordinary cusps. Since only 
one type of singularities is present  and (\ref{quasiadjcusp}) shows that 
there is only one constant quasiadjunction
$\kappa={1 \over 6}$, the Alexander polynomial 
has the form 
$$[(t-exp({{2 \pi i} \over 6}))(t-exp(-{{2 \pi i} \over 6}))]^s=
(t^2-t+1)^s$$ 
Now the linear system in question consists on the 
curves having degree $6-3-{6 \over 6}=2$ with local equations 
belonging to the maximal ideals of the singular points. 
Since the dimension of the space of quadrics in 6, the expected 
dimension of our linear system is 0 and if a quadric containg
all six cusps does exist then the 
actual dimension is 1 (one can show that this is the maximal possible value).
Hence $s=1$ and the Alexander polynomial is $t^2-t+1$.

For a sextic $\hat C_3$ with nine cusps, 
which is dual to a non singular cubic, 
one has ${{\rm dim} H^0(\PP^2, {\cal J}_{\kappa_i}
(d-3-d\kappa_i))}=0$ and $\chi({\cal J}_{\kappa_i})=-3$
and hence:

\begin{equation} 
\Delta_{\hat C_3}(t)=(t^2-t+1)^3
\end{equation}

For the curve from $\cite{turpin}$ given by the equation:
$f^2_{3n}+f^3_{2n}=0$ where $f_n$ is a generic form of degree $n$, which 
has only ordinary cusps at $6n^2$ points forming a complete 
intersection of curves of degrees $2n$ and $3n$ the exponent
of $t^2-t+1$  in theorem \ref{positioncurves} is the superabundance
of the curve of degree $6n-3-{{6n} \over 6}$ passing through 
this complete intersection. By a theorem of Cayley-Bacharach
this superabundance is 1 and hence the Alexander polynomial 
is $t^2-t+1$.

\subsection{Commutative fundamental groups}\label{commfund}

\subsubsection{Commutativity in terms of local type of 
singularities. Nori's theorem.} 

Historically, much of the work on the fundamental groups 
of the complements, was focused on the cases when 
the fundamental group is abelian. In this case Prop. \ref{homology}
yields the complete calculation of $\pi_1$. For example, as was pointed
out in the introduction, 
F.Severi was claiming
that the fundamental group of the complement to a curve having 
nodes as the only singularities is abelian. More precisely he claimed 
the irreducibility of the stratum of nodal curves (this was proven by 
much later by 
J.Harris in \cite{harris}). The irreducibility of this stratum yields that 
each nodal curve can be degenerated into a union of lines in general 
position and for such union (these days called a generic 
arrangement of lines) a direct calculation shows that 
the fundamental group of the complement is free abelian.
More generally than in the case of nodal curves,
one expects, speaking very vaguely,  
that if a curve has not too many singularities or 
if the singularities are sufficiently mild then the fundamental groups 
of the complement will be abelian. 
A precise result in this direction follow from a theorem of M.Nori:

\begin{theo} Let $D$ and $E$ be a curves on a non singular 
surface $X$. Assume that $D$ has nodes as the only singularities,  
that $D$ and $E$ intersects transversally and that for an irreducible 
component $C$ of $D$ one has $C^2>2r(C)$ where $r(C)$ is the number 
of nodes on $C$. Then 
$N={\rm Ker} \pi_1(X-D-E) \rightarrow \pi_1(X-E)$ is abelian.
\end{theo}

For plane curves one obtains the following which   
extends earlier commutativity results of S.Abhyankar.

\begin{theo} 
For a germ $\phi$ of a curve singularity in $\CC^2$
let 
us define the invariant $e(\phi)$ as follows. Let 
$\Phi: S \rightarrow \CC^2$ be a resolution of the singularity
and $\Phi^*(\phi)=F+G$ where $F$ is the proper transform of 
$\phi=0$, $G$ is the exceptional set and $F$ and $G_{red}$ meet transversally.
Let $e(\phi)=G(G+2F)$ and let $F(C)$, for a curve $C$ on a non singular 
projective surface, be the sum of invariants $e(\phi)$ 
for all singularities of $C$. 
If $C^2 > F(C)$ then the extension $\pi_1(X-C) \rightarrow \pi_1(X)$ 
is central. 
\end{theo}

{\it Proof.} Apply Nori's theorem the proper transform $C'$ of a resolution 
of singularities of $C$. Then if $G$ is the exceptional set 
then $C^2=(C'+G)^2={C'}^2+2(C',G)+G^2={C'}^2+F(C)$.
Hence the assumed inequality translates into ${C'}^2>0$. Hence Nori's theorem 
yields the conclusion. Note that for node we have $G=2E$ where $E$ is the 
exceptional line and $C'=L_1+L_2$. Hence 
$G^2+2(G,C')=4E^2+2\cdot 2E(L_1+L+2)=4$. For a cusp $F(\phi)=6$. In particular
on a simply-connected surface the fundamental group of the complement
to a curve with $\delta$ nodese and $\kappa$ cusps
is abelian if $C^2>6 \kappa+4\delta$.

The following question is still open:
\begin{quest} Let $N$ be a normal subgroup of $\pi_1(X)$ 
generated by the images of the fundamental groups of non 
singular models of components. Does $N$ has a finite index in $\pi_1(X)$
\end{quest}
If so, then the fundamental group of a surface, containing 
a rational curve with positive self-intersection, must be finite.

\subsubsection{On a  proof of Nori's theorem.}

Let us consider a special case when  $E=\emptyset$, and 
$C$ is an irreducible nonsingular curve on $X$. 
Let $U$ be a tubular neighborhood 
of $C$. Then $U-C \rightarrow C$ is a circle fibration and the fiber $\delta$
 is the element of $\pi_1(U-C)$ belonging to the center of the latter
group.
Since in this case the assumption of the theorem is $C^2>0$, the
 theorem of Nakai and Moishezon (cf. \cite{hartshorne}) 
yields that $C$ is ample and hence 
a small deformation $D$ of $nC$, which we may assume belongs to 
$U$, is very ample and also smooth. By usual Lefschetz theorem 
$\pi_1(D-C) \rightarrow \pi_1(X-C)$ is surjective and hence $\pi_1(U-C) 
\rightarrow \pi_1(X-C)$ is surjective as well. Therefore the image if the 
class of $\gamma$ in $\pi_1(X-C)$ belongs to its center.
On the other hand, any element in 
$N={\rm Ker} \pi_1(X-D) \rightarrow \pi_1(X)$ 
is product of elements conjugated to $\gamma$. Indeed, take
such element $\delta$ and consider 2-disk $\Delta$ which it bounds in $X$.  
We can assume that $\Delta \cap C$ consists of finitely many transversal 
intersections. Therefore $\delta=\Pi \delta_i$ where 
$\delta_i=\alpha_i \gamma_i \alpha^{-1}$ with $\gamma_i$ being 
a fiber of $U-C \rightarrow C$ and $\alpha_i$ is a path going 
from the base point to a point on the boundary of $U$. In particular  
$\delta_i$ is conjugate to $\gamma$ in $\pi_1(X-C)$ and hence 
is equal to $\gamma$. Hence $\delta$ is a power of $\gamma$ i.e. 
$N$ is cyclic.

Crutial in the proof of Nori's theorem in 
the case of nodal $C$ is the following Nori's weak Lefschetz theorem
which is very interesting by itself.

\begin{theo} Let $i: H \rightarrow U$ be an embedding of 
 connected compact complex analytic
subspace (possibly non reduced) into a connected complex manifold $U$ 
in which $H$ is defined by a locally principal sheaf of ideals. Assume that 
$\O_U(H) \vert H$ is ample and that $dim U >2$. 
Let $q: U \rightarrow X$ be holomorphic local isomorphism with the target
being a smooth projective variety and $h=q \circ i$. 
Let $R$ be an arbitrary Zariski 
closed subset and $G={\rm Im} \pi_1(U-q^{-1}(R)) \rightarrow \pi_1(X-R)$.Then 
$G$ is a subgroup of finite index.  
\end{theo}

\subsection{Higher homotopy groups}  

Another natural invariants of the homotopy type of the 
complement are the higher homotopy groups of the complement.
However for curves, the higher homotopy groups, unlike the
fundamental groups, it seems, 
do not have an algebro-geometric significance. 
Moreover, in most cases the higher homotopy groups, 
considered as abelian groups are infinitely generated. 
A more useful way to consider them is by using the action of 
$\pi_1$ on $\pi_k$ i.e. consider $\pi_k$ as a module 
over $\pi_1$. However unless $\pi_1$ is abelian, 
understanding modules over $\pi_1$ 
involves a subtle non commutative algebra. 
For curves however, as will be explained in the next section, 
the homotopy type of the complement
is determined by another invariant of the pair $(\PP^2,C)$ i.e. 
the {\it braid monodromy}.
On the other hand for hypersurfaces in $\PP^{n+1}$ with $n>1$ the 
homotopy groups in dimensions up to $n$ have interesting
algebro-geometric meaning which we shall 
proceed to discuss.

\subsubsection{Action of the fundamental group on higher homotopy 
groups} 

Let us start with the example which shows why the homotopy 
groups of simplest topological spaces are infinitely generated.

\begin{exam} Let us consider $\pi_2(S^1 \vee S^2)$. Clearly 
$\pi_1(S^1 \vee S^2)=\ZZ$. On the other hand $\pi_2(S^1 \vee S^2)$
can be identified with $\pi_2$ of the universal cover of 
$S^1 \vee S^2$. Viewing the universal covering map of the 
circle as the the quotient of $\RR$ by the subgroup of integers 
makes it natural to view the universal cover of $S^1 \vee S^2$ 
as the real line with $S^2$'s attached at the integer points.
Hence the universal cover has $H_2$, and by Hurewicz theorem $\pi_2$,
isomorphic to $\ZZ^{\infty}$. On the other hand, 
since the deck transformation of the universal cover 
acts transitively on $S^2$'s attached to $\RR$, both $H_2$ and $\pi_2$ 
are cyclic modules over the group of deck transformations i.e 
$\pi_2(\widetilde {S^1 \vee S^2})=\ZZ[t,t^{-1}]$ ($\widetilde X $ 
denotes the universal cover).
\end{exam}

In general, the homotopy groups can be given the structure of a 
module over the fundamental group using the Whitehead product:
$\pi_n \times \pi_m \rightarrow \pi_{n+m-1}$.  
In the case when $\pi_i(X)=0$ for 
$2 \le i \le n-1$, if $\tilde X$ is the universal 
cover then $\pi_n(X)=\pi_n(\tilde X)=H_n(\tilde X)$
and the action of $\pi_1(X)$ is just the 
action of the deck transformations on the homology.

Such $X$ come up naturally:

\begin{theo}\label{lowhomotopyvanish} Let $V$ be a hypersurface in 
$\PP^{n+1}$ having only isolated singularities.
Let $H$ be a generic hyperplane. 
Then $\pi_1(\PP^{n+1}-V \cap H)=\ZZ$  and 
$\pi_i(\PP^{n+1}-V \cap H)=0$ for $2 \le i \le n-1$.
Moreover, $\pi_n(\PP^{n+1}-V \cap H) \otimes \QQ$ is 
a $\QQ[t,t^{-1}]$-torsion module. 
\end{theo}

More generally, the Lefschetz hyperplane section theorem 
yields that the conclusion of the theorem holds 
for arbitrary hypersurfaces in $\PP^N$ for which the 
singular locus has codimension $n+1$. To see this (and also the 
first part of theorem \ref{lowhomotopyvanish}) recall 
it:

\begin{theo}({\bf Lefschetz hyperplane section theorem})

(a) Let $X$ be a projective subvariety having dimension $n$ and 
let $L$ be a codimension $d$ linear subspace such that
$X$ is a local complete intersection outside of $L$.
Then $$\pi_i(X \cap L) \rightarrow \pi_i(X)$$
is isomorphism for $0 \le i < n-d$ and surjective for 
$i=n-d$.

(b)Let $X$ be a quasiprojective. The conclusion of (a) 
take place for generic $L$.  

\end{theo}

Vanishing statement in theorem \ref{lowhomotopyvanish} follows 
from this and calculation of the homotopy groups of the 
complement to non singular hypersurfaces.

Recently, L.Maxim (\cite{maxim})
showed that the homology of infinite cyclic covers of the 
complement to an affine hypersurface, generic relative 
to hyperplane at infinity, are torsion modules in all dimensions 
except the top one.

\subsubsection{Orders of the homotopy groups}

It follows from the theorem  \ref{lowhomotopyvanish} and the 
classification of modules over PIDs that 
$$\pi_n(\PP^{n+1}-V \cap H) \otimes \QQ=\oplus 
\QQ[t,t^{-1}]/\Delta_i(t)$$ for some polynomials $\Delta_i(t)$
defines up to a unit in $\QQ[t,t^{-1}]$. We call $\Delta(t)=\Pi_i \Delta_i(t)$
the order of the group $\pi_n$. Though $\Delta(t)$ cannot 
be calculated in terms of a local data of singularities there is the
following divisibility relation, which generalizes the
divisibility relation for the Alexander polynomials:

\begin{theo}\label{divisone}
({\bf Divisibility theorem I})
The order of $\pi_n(\CC^{n+1}-V)$ divides the 
product of characteristic polynomials of the 
monodromy operators of singularities of $V$:
$$ \Delta(\CC^{n+1}-V) \vert \Pi_{P_i \in Sing(V)}
 \Delta_{P_i}(t)$$
\end{theo}

Note that as it stated, one should assume that
$V$ it transversal to the hyperplane at infinity. 
However one can define correction factors corresponding 
to the singularities at infinity so that, 
if one multiplies by these correction factors the right 
side in \ref{divisone}, the divisibility relation
will hold. 

\begin{theo}({\bf Divisibility theorem II}) Let $V$ be a hypersurface
transversal to the hyperplane at infinity $H_{\infty}$. Let $S_{\infty}$ 
be the boundary of a small tubular neighborhood of $H_{\infty}$ and 
let $L_{\infty}=V \cap S_{\infty}$. Then the homology of the infinite 
cyclic cover of $S_{\infty}-L_{\infty}$ is a torsion $\CC[t,t^{-1}]$-module
and $\Delta_{\infty}$ and $\Delta (\CC^{n+1}-V) \vert  \Delta_{\infty}$.
\end{theo}
(see \cite{annals} for a statement in the case with a 
weaker than transversality to $H_{\infty}$
assumption)

\section{Homotopy groups via pencils.}\label{homotopygroupsviapencils}

\subsection{Van Kampen theorem and braid monodromy}

Now let us consider how one actualy can   
calculate the fundamental 
group of a complement in the case of curves and how 
to calculate the first non trivial homotopy group of the 
complement in the case of hypersurfaces. 
In this section we shall deal with the curves
(cf. also \cite{eyral} where the case of possibly singular 
quasiprojective varieties is discussed).

Let $C$ be a curve on a projective surface $X$ for which 
we want to describe $\pi_1(X-C)$. Consider 
a line bundle $\L$ on $X$ such that 
${\rm dim} H^0(X,\L)) \ge 2$ and select a 2-dimensional
linear system ${\bf L} \subseteq  H^0(X,\L)$. Let $B$ be the 
base locus of $\bf L$ (it contains at most $c_1(\L)^2$ points).
We shall assume for simplicity that $B \cap C =\emptyset$.
The classical case is $X=\PP^2$, $\L=\O(1)$ and 
${\bf L} \subset H^0(\PP^2,\O(1))$ consists of sections with 
the zerosets containing  
a fixed point.  We have a regular map onto the projectivization of 
$\bf L$:  
\begin{equation}\label{pencil}
p: \ X-B \rightarrow \PP({\bf L})=\PP^1
\end{equation}
with generic fiber $L_{t_0}-L_{t_0} \cap B, t_0 \in \PP^1$ 
being non singular by Bertini's (or Sard's) theorem. 
Though generic element of ${\bf L}$ may be singular at points of $B$,
we shall make additional assumtion that $L_{t_0}$ is 
non singular at any $p \in B \cap C$.
The curve $L_{t_0}$ is ample and 
hence $\pi_1(L_{t_0}-L_{t_0} \cap C) \rightarrow
\pi_1(X-C)$ is surjective by Lefschetz theorem. We want to
describe the kernel of this map. Let $Sing \subset \PP^1$ be the 
(finite) subset of points $t_1,...,t_N$ 
corresponding to singular members of 
the pencil. Each fiber of the pencil (\ref{pencil}) 
is a punctured curve (which, if $L_{t_0}$ is non singular at 
the points of $B$, has genus 
$g(L_{t_0})={{c_1(\L)(K_X+c_1(\L))} \over 2}+1$).

For each $d$ one can define the braid group $B_d(L_{t_0}-B)$ which 
is the group of isotopy classes of orientation 
preserving diffeomorphisms of  $L_{t_0}$ 
which are constant in a neighborhood of $B$ in $L_{t_0}$.
In the case $L_{t_0}-B \cap L_{t_0}=\CC$ one obtains the classical 
Artin's braid group with generators $\sigma_i, i=1,...,{d-1}$
and relations
\begin{equation} 
\sigma_i \sigma_j=\sigma_j \sigma_i \ \vert i-j \vert \ge 2, 
\sigma_i \sigma_{i+1} \sigma_i=\sigma_{i+1} \sigma_{i} \sigma_{i+1}\ \ 
i=1,..,d-2
\end{equation}
(for presentations of braid groups similar to $B_d(L_{t_0}-B)$ 
by generators and relations and 
extending this 
one, see \cite{scott})

We want to construct a homomorphism $\pi_1(\PP^1-Sing, t_0) \rightarrow
B_d(L_{t_0}-B)$ called {\it braid monodromy} and use it to describe
$\pi_1(X-C)$. We shall start by defining ``good'' systems of generators of 
$\pi_1(\CC-Sing)$, which we shall use to give 
a finite presentation for this fundamental group.

\begin{dfn}\label{goodcollection}
Let $Sing=\{t_1,...,t_N\}$.
A system of generators $\gamma_i \in \pi_1 ({\bf C}-
\bigcup _i t_i, t_0)$ is called {\it good} if each of the loops
 $\gamma_i: S^1
\rightarrow  {\bf C}- \bigcup_i t_i$ extends to a map of the disk
$D^2 \rightarrow  {\bf C}$ with non-intersecting images for distinct
$i$'s.
\end{dfn}
 
One way to construct a good system of generators is
the follwoing. Select a system of small disks $\Delta _i$ about each 
point $t_i$
$i=1,...,N$, and
choose a system of $N$ non-intersecting paths $\delta_i$ connecting the
base point $t_0$ with a point of $\partial \Delta _i$.  Then
$\gamma_i = \delta ^{-1} \circ \partial \Delta _i \circ \delta_i$
is a good system of generators
(with, say, the counterclockwise orientation of $\partial \Delta _i$).
We shall need also good systems of generators of the fundamental groups
of the complements to a finite set of $N$ points on a {\it compact} Riemann 
surface having genus $g \ge 0$. 
Those are the systems of generators $\gamma_1',...,
\gamma_{2g}',$ consisting of the 
images $2g$ sides of a $4g$-gon for some presentation of the surface 
as a $4g$-gon with identified sides and a good system of generators 
$\gamma_1,...,\gamma_N$ of 
the complement to $N$ points in this $4g$-gon in the above sense.
We have the only relation 
\begin{equation}\label{grouprelation}
R: \ \Pi \gamma_1 \cdot ...\cdot \gamma_N=
\Pi[\gamma_i',\gamma_{i+1}']
\end{equation}
 In the case $g=0$ this relation 
becomes $\Pi \gamma_1 \cdot ...\cdot \gamma_N=1$.

Now let us define the braid monodromy corresponding to an element 
$\gamma \in \pi_1(\PP^1-{\rm Sing})$.
Let $\gamma \in \PP^1-{\rm Sing} $ be the image of an embedding of $S^1$ 
taking the base point to $t_0$. We can view $\gamma$ as the image of the 
map $\iota: I \rightarrow \PP^1-Sing$ ($I$ is the unit interval)
such that $\iota(0)=\iota(1)=t_0$. Then $(X-B-C) \times_{\PP^1-Sing} I$ 
is a locally trivial fibration over $I$ and hence is a trivial fibration. 
This means
that there is a map $\Phi: L_{t_0}-L_{t_0} \cap C 
\times I \rightarrow X-B$ such that 
$\Phi(t) \vert_{L_{t_0}-L_{t_0} \cap C \times t}$ 
is a homeomorphism onto $L_{t}-L_{t} \cap C$. 
Note that though 
$\Phi$ is not unique any two choices are isotopic via isotopies 
commuting with projections on $I$. Hence we obtain the map
$\Phi (1): L_{t_0}-L_{t_0} \cap C \rightarrow L_{t_0}-L_{t_0} \cap C$
and the isotopy class of this map is well defined. 
We can assume that this map keeps $B$ fixed.
One checks immediately that dependence on $\iota$ yields homotopic
maps $\Phi(1)$ and a homotopy of $\gamma$ extends to a homotopy of $\Phi(1)$
(but $B$ may not be possible to preserve). 
Hence we obtain the braid monodromy homomorphism:
\begin{equation}\label{monodromyhomo} 
\pi_1(\PP^1-Sing) \rightarrow \pi_0(Diff(L_{t_0}))=B_d(L_{t_0})
\end{equation}
(where $d=(C,L_t)$ and the last group is the braid group of Riemann surface
$L_{t_0}$. 

There is a useful way to encode algebraically the homomorphism
(\ref{monodromyhomo}) using the choice of a 
good system of generators of $\pi_1(\PP^1-{\rm Sing})$. 
Let us fix a fiber $L_{t_{\infty}}$
of the pencil which we shall call the fiber at infinity. Then we can select
monodromy transformations all fixing a neighborhood of $B$ for all $\gamma_i$ 
i.e. we obtain ordered system of braids: 
$\beta(\gamma_i)=\Phi_{\gamma_i}(1) \in B_d(L_{t_0}-L_{t_0} \cap C-B)$
with the order given by the order of the good systems of generators.
The latter is given by the counterclickwise ordering of loops 
about the point $t_0$. Moreover, the product is a fixed word
in $B_d(L_{t_0})$ independend of $C$.
For example we obtain in the case of curves in $\CC^2$:
\begin{equation}  
\Pi\beta(\gamma_i)=\Delta^2   
\end{equation}
where $\Delta^2$ is the generators of the center of the 
Artin's braid group $B_d$ (cf. \cite{bowden}).

We have the following calculation in terms of the braid monodromy
originated by Zariski-van Kampen:

\begin{theo}
Let $b \in \partial T(B) \cap L_{t_0}$ where $T(B)$ is a neighborhood
of $B$ in $X$ and let $\alpha_j$ be 
a good system of generators of $\pi_1(L_{t_0}-L_{t_0} \cap C,b)$.
Let $R$ be the relation among $\alpha_j$. Then
$$ \pi_1(X-C-L_{t_{\infty}})=
\pi_1(L_{t_0}-L_{t_0}-B,b)/(\beta(\gamma_i)(\alpha_i)\alpha_i^{-1})$$
(quotient by the normal subgroup generated by specified elements).
The group $\pi_1(X-C)$ can be obtained by adding to the above 
the relation $R$.
\end{theo}

In the case of plane curves we have just the homomorphism into 
Artin's braid group which by itself is an interesting invariant
of plane curves
containing more information than the fundamental group.
For example the braid monodromy 
determines the homotopy type of the complement $\CC^2-C$ 
(cf. \cite{homotopytype}).
Many calculations are done for 
curves $C$ which are the branching curves of generic projections of 
surfaces (cf. \cite{teicher}).
Recently braid monodromy found applications in
 symplectic geometry (cf. \cite{symplectic}).

\subsection{Homotopy groups via pencils}

Now let $V$ be a hypersurface in $\CC^{n+1}$ transversal to 
the hyperplane at infinity and having only isolated singularities. 
We want to describe calculation of the first non trivial 
homotopy group $\pi_n(\CC^{n+1}-V)$ in terms of pencils generalizing 
the Zariski-van Kampen procedure described above.

We start with a high dimensional analog of the braid group 
and a linear representation generalizing the Burau representation 
of the braid group. In higher dimensions we have several 
candidates for such a generalization.

Let us consider a sphere $S^{2n-1}$ in ${\bf C}^{n}$ 
of a sufficiently large radius. Let
$\partial_ \infty V= V \cap S^{2n-1}$ and let  
$Emb (V,{\bf C}^{n})$ be the space of submanifolds of ${\bf C}^n$ 
with the following property: each is 
diffeomorphic to $V$ and, moreover, is isotopic to the chosen embedding
 of $V$. In addition we require that
 for any $V^{\prime} \in Emb (V,{\bf C}^{n})$ one has $V^{\prime}
  (V) \cup S^{2n-1}=\partial _\infty V$.
We shall use the topology with the basis consisting of 
sets $U(V,\epsilon)$  of submanifolds $V' \subset \CC^n$ 
which belong to the tubular neighborhood of $V$ having radius $\epsilon$
and which are isotopic to $V$. 

Let us  describe a linear
 representation 
\begin{equation}\label{firstmap}
\pi_1 (Emb (V,{\bf C}^{n})) \rightarrow Aut 
\pi_n ({\bf C}^{n}-V)
\end{equation}
 After a choice 
of a basis in the $\pi_1 ({\bf C}^{n}-V)$-module $\pi_n ({\bf C}^{n}-V)$
this homomorphism becomes the 
 homomorphism into $GL_r ({\bf Z} [t,t^{-1}])$ where
 $r$ is the rank of $\tilde H_n ({\bf C}^{n}-V,{\bf Z})$ (the reduced
 homology of the complement). It is given in terms of the representation of 
another group which also is a candidate for the high-dimensional braid group.

Let $Diff ({\bf C}^{n},S^{2n-1})$ be the
group of diffeomorphisms of ${\bf C}^{n}$ acting as the identity
 outside $S^{2n-1}$. This group can be identified with the group 
$Diff (S^{2n},D^{2n})$ of the diffeomorphisms of the sphere fixing
 a disk. Let $Diff ({\bf C}^{n},V)$ be the subgroup
 of $Diff ({\bf C}^{n},S^{2n-1})$ of the diffeomorphisms  which
 take $V$ into itself. 

The group $Diff ({\bf C}^{n},S^{2n-1})$ acts
 transitively on $Emb (V,{\bf C}^{n})$ with the stabilizer $Diff ({\bf
 C}^{n},V)$.  Therefore we have the following exact sequence:
 \begin{equation}\label{longsequence} 
\pi_1 (Diff ( S^{2n+2},D^{2n+2})) \rightarrow  \pi_1 (Emb ({\bf C}
 ^{n},V)) \rightarrow \pi_0 (Diff ({\bf C}^{n},V)) 
\end{equation}
$$\rightarrow
\pi_0  (Diff( S^{2n+2},D^{2n+2})) \rightarrow $$
 Any element in $Diff ({\bf C}^{n},V)$ induces the self map  of
 ${\bf C}^{n}-V$ and also the self map of the universal (in the case
  $n=1$ universal cyclic)
  cover of this space. Hence it induces an automorphism of $H_n
(\widetilde {{\bf  C}^n-V,{\bf Z}})=\pi_n ({\bf C}^n-V)$, $n>1$.
This gives the representation:
\begin{equation}\label{secondrep}
\pi_0 (Diff ({\bf C}^{n},V)) \rightarrow  Aut 
\pi_n ({\bf C}^{n}-V)
\end{equation}

 The composition of the boundary
 homomorphism in (\ref{longsequence})
 with the map (\ref{secondrep})
results in representation (\ref{firstmap}).
The groups $\pi_1 (Emb ({\bf C}
 ^{n},V))$ and $\pi_0 (Diff ({\bf C}^{n},V))$
are high-dimensional analogs of the braid groups and their 
algebraic study was not carried out so far. However 
some high-dimensional analogs of the mapping class groups 
were studied in (cf. \cite{krylov}).

\par  In the case $n=1$, $V$ is just a collection of points in ${\bf
C}$,
 $\pi_1 (Emb ({\bf C},V))=\pi_0 (Diff ({\bf C},V))$ is  Artin's
 braid group, and this construction gives the homomorphism of the braid
 group into $Aut H_1 (\widetilde {{\bf C}-V,{\bf Z}})$ which, after a
choice  of the basis in $ H_1 (\widetilde {{\bf C}-V})$
 corresponding to the choice
 of the generators of the braid group, gives the reduced Burau
representation. In higher dimensions the isomorphism
$\pi_1 (Emb (V,{\bf C}^{n}))= \pi_0 (Diff ({\bf C}^{n},V))$ fails.

\begin{problem} Calculate the groups $\pi_1 (Emb (V,{\bf C}^{n}))$
and  $\pi_0 (Diff ({\bf C}^{n},V))$
\end{problem}

\par Now we  can define the relevant monodromy operator corresponding to
a  loop in the parameter space of a linear pencil of hyperplane sections.
  By our assumptions, the projective closure of $V$ is 
a hypersurface in ${\bf P}^{n+1}$
 which has only isolated singularities. Let $H$ be the hyperplane at
 infinity (which is transversal to $V$).
  Let $L_t$, $t \in {\bf C}$, be
 a pencil of hyperplanes the projective closure of which has as the base
locus a hyperplane $M$ in $H$ such that $M$ also is transversal to $V$. Let $t_1,..
.,t_N$ denote those $t$ for which $V \cap L_t$ has a singularity.
We shall assume that for any $i$ the singularity of $V \cap L_{t_i}$ is
outside of $H$. The pencil   $L_t$
over ${\bf C}- \bigcup_i t_i$  defines a 
locally trivial fibration $\tau$ of ${\bf C}^{n+1}-V$
 with a non-singular hypersurface in ${\bf C}^n$ transversal
to the hyperplane at infinity as a fiber. The restriction of this  fibration
 on the complement to a sufficiently large ball is trivial, as follows
 from the assumptions on the singularities at infinity.
Let $\gamma: [0,1] \rightarrow  {\bf C}- \bigcup_i t_i$ $(i=1,...,N)$
be a loop with the base point $t_0$.
   A choice of a trivialization of the pull back  of the
fibration $\tau$ on $[0,1]$ using $\gamma$, defines a loop $e_{\gamma}$
in $Emb (L_{t_0},V \cap L_{t_0})$. Different trivializations
 produce homotopic loops in this space.
\begin{dfn} The monodromy operator corresponding
 $\gamma$ is the  element in 
\par \noindent $Aut (\pi_n (L_{t_0}-L_{t_0} \cap V))$
 corresponding in (\ref{firstmap}) to  $e_\gamma$
\end{dfn}
\par Next we will need to associate the following homomorphism with a singular
fiber $L_{t_i}$
 and a loop $\gamma$ with the  base point $t_0$
 in the parameter space of the pencil where $\gamma$ bounds 
 a disk $\Delta_{t_i}$ not containing other singular points
  of the pencil :
\begin{equation}\label{degoperator}
   \pi_{n-1} (L_{t_i}-L_{t_i} \cap V) \rightarrow
 \pi _n (L_{t_0}-L_{t_0} \cap V)/ Im (\Gamma -I). 
\end{equation}
 Here $\Gamma$ is the monodromy operator corresponding to $\gamma$.
 \par To construct (\ref{degoperator})
 let us note that the $\pi_1$-module on the right in (\ref{degoperator}) is
 isomorphic to the homology $H_n (\widetilde
  {\tau ^{-1} (\partial \Delta _{t_i})} ,{\bf Z})$ of the infinite
 cyclic cover of the restriction of the fibration $\tau$ on the
 boundary of $\Delta _{t_i}$. This follows  immediately  from the Wang
sequence
 of a fibration over a circle and the vanishing of the homotopy of
 $L_{t_0}-L_{t_0} \cap V$ in dimensions below $n$. Let $B_i$ be a
 polydisk in ${\bf C}^{n+1}$ such that $B_i=\Delta _i \times B$ for
 a certain polydisk $B$ in $L_{t_0}$. Then $\widetilde {\tau ^{-1}
  (\Delta _i)-B_i}$ is a trivial fibration over $\Delta_i$ with the
 infinite cyclic cover $\widetilde {L_{t_i}-L_{t_i} \cap V}$ as a fiber.
 In particular,  one obtains the map:
  $$\pi_ {n-1} (L_{t_0}-L_{t_0} \cap V)=
 H_{n-1} (\widetilde {L_{t_0}-L_{t_0} \cap V},{\bf Z}) \rightarrow
 H_n (\widetilde {\tau ^{-1} (\Delta _i)-B_i},{\bf Z})=$$
\begin{equation}\label{degoperator2}
H_{n-1} (\widetilde{L_{t_0}-L_{t_0} \cap V},{\bf Z}) \oplus
 H_n (\widetilde {L_{t_0}-L_{t_0} \cap V},{\bf Z}) 
\end{equation}
\begin{dfn}\label{defdeg}
 The degeneration operator is the map (\ref{degoperator})
 given by composition of the map (\ref{degoperator2}) with  the map
 $H_n (\widetilde {\tau ^{-1}
   (\Delta _i)-B_i},{\bf Z}) \rightarrow H_n(\widetilde {\tau^{-1}
 (\Delta _i)},{\bf Z})=\pi _n (L_{t_0}-L_{t_0} \cap V)$ induced by
 inclusion.
\end{dfn}

The following is a high-dimensional analog of the van Kampen theorem.

\begin{theo}Let $V$ be a hypersurface in ${\bf P}^{n+1}$
having only isolated singularities and transversal to the hyperplane $H$
at infinity. Consider a pencil of hyperplanes in ${\bf P}^{n+1}$
the base locus $M$ of which belongs to $H$ and is transversal in $H$ to
$V \cap H$. Let
${\bf C}^n_t$ ($t \in \Omega$) be the  pencil of hyperplanes in ${\bf
C}^{n+1}={\bf
P}^{n+1}-H$ defined by  $L_t$ (where $\Omega={\bf C}$ is the set
parameterizing all elements of the pencil $L_t$ excluding $H$). Denote by
$t_1,...,t_N$  the collection of
those $t$ for which $V \cap L_t$ has a singularity.
We shall assume that the pencil was chosen so that $L_t \cap H$ has at
 most one singular
  point outside of $H$. Let $t_0$ be different from either of $t_i$
($i=1,...,N$). Let
$\gamma_i$ ($i=1,...,N$) be a good collection, in the sense described
in definition (\ref{goodcollection}), 
of  paths in $\Omega$ based in $t_0$ and forming a
basis of
$\pi_1 (\Omega- \bigcup_i t_i,t_0)$ and let
$\Gamma _i \in Aut (\pi_n ({\bf C}^n_t-V \cap {\bf C}^n_t))$   be the
monodromy automorphism corresponding to $\gamma_i$. Let $\sigma_i:
\pi_{n-1} ({\bf C}^n_{t_i}-V \cap {\bf C}^n_{t_i}) \rightarrow
\pi_n ({\bf C}^n_{t_0} - V \cap {\bf C}^n_{t_0}) ^{\Gamma_i}$ be the
 degeneration operator of the homotopy group of a special
 element of the pencil
into the corresponding quotient of covariants constructed above. Then
 $$\pi_n ({\bf C}^{n+1}-V \cap {\bf C}^{n+1})= \pi_ n ({\bf C}^n-V \cap
 {\bf C}^n) / (Im (\Gamma_1-I),Im \sigma_1,...,
 Im (\Gamma_N -I), \sigma_N )$$
\end{theo}

We refer for a proof to the paper \cite{annals}.
There is another way to describe this homotopy group 
replacing the degeneration operator by a variation 
operators on the homotopy groups which we shall describe now.

\subsection{Variation operators}

Variation operators classically defined in homology (or cohomology).
The idea of defining homotopy variation operator 
comes from the fact that the homotopy groups on question
are the homology groups of covering spaces. 
A description of the homotopy groups using variation operators 
was carried out in \cite{denisme}.

We shall continue to use the notations from previous section but in 
addition let us select 
$e \in M-M \cap V$ which we shall use as the base point 
for the homotopy groups. 
The homotopy variation operator is a certain homomorphism
of $\ZZ[t,t^{-1}]$-modules:
\begin{equation}
 {\cal V}_i \colon \ \pi_n(L_{t_0}-L_{t_0}\cap V,M-M\cap V,e)
 \longrightarrow \pi_n(L_{t_0}-L_{t_0}\cap V,e)
\end{equation}
associated with each $\gamma_i$ for $1\leq i\leq N$.
 
As in definition (\ref{defdeg}) of degeneration operators 
 we shall go to the $d$ fold cover
 and use the homological variation operators.
Let $W \subset \PP^{n+2}$ 
be the $d$-fold cyclic branched over $V$ cover of $\PP^{n+1}$, 
$j: V \rightarrow W$ the embedding and $\L_t$ be the pull back of 
the pencil $L_t$ on $W$. By abuse of notations we shall denote
by the same letter the hyperplane in $\PP^{n+2}$ cutting the corresponding 
divisor on $W$. 
$\L_{t_0} \cap W$ is the $d$-fold cover of $L_{t_0}$ branched over 
$V \cap L_{t_0}$. We shall consider the homological 
variation operators:
\begin{equation}
 V_i\colon \ H_n(\L_{t_0} \cap(W-j(V)), {\cal M} \cap(W-j(V))\bigr)
 \longrightarrow H_n(\L_{t_0}\cap(W-j(V)))
\end{equation}
associated, for $1\leq i\leq N$, with the homotopy classes $\gamma_i$.

The definition and the properties of operators $V_i$ 
are discussed in \cite{denisLondon}.
For a relative $n$-cycle $\Xi$ on $\L_{t_0}\cap(W-j(V))$ with boundary in
$\M \cap(W-j(V))$, one defines ($[ \cdot ]$ denotes the class of a cycle):
\begin{equation}
\label{Vardef}
 V_i([{\Xi}]_{(\L_{t_0}\cap(W-j(V)),\M \cap(W-j(V)))})
 =[H_{i*}(\Xi)-\Xi]_{\L_{t_0}\cap(W-j(V))}
\end{equation}
where $H_i$ is the geometric monodromy corresponding to $\gamma_i$.
Since  $H_i$ leaves the points of
$\M \cap(W-j(V))$ fixed the chain
$H_{i*}(\Xi)-\Xi$ is actually an absolute cycle and the
correspondence $\Xi\mapsto H_{i*}(\Xi)-\Xi$ induces a
homomorphism $V_i$ at the homology level (\cite[Lemmas 4.6
and 4.8]{denisLondon}). 
This homomorphism depends only on the homotopy
class $\gamma_i$ (\cite[Lemma~4.8]{denisLondon}).

Now, if $n\geq 2$ then for $1\leq i\leq N$,
using the isomorphism $\alpha_{t_0}$ and the
homomorphism $\bar\alpha_{t_0}$, 
$V_i$ yields the homotopical variation operator ${\cal V}_i$ by
requiring that the following diagram will be commutative:
\begin{equation}
\label{Valpha}
\matrix{H_n (\L_{t_0}\cap(W -j(V)),\M \cap(W-j(V))) & \buildrel V_i
\over \longrightarrow & H_n (\L_{t_0}\cap(W-j(V))) \cr 
\uparrow \bar \alpha_{t_0} & & \uparrow \alpha_{t_0} \cr 
 \pi_n(L_{t_0}-L_{t_0}\cap V,M-M\cap V,e) & \buildrel {\cal V}_i \over
 \longrightarrow & \pi_n (L_{t_0}-L_{t_0}\cap V,e). \cr} 
\end{equation}
As $V_i$ depends only on the homotopy class $\gamma_i$ 
so do the operators ${\cal V}_i$.

With these definitions we have the following (cf. \cite{denisme}):

\begin{theo}
Let $V$ be a hypersurface in $\PP^{n+1}$ with $n\geq2$ having only
isolated singularities. Consider a pencil $(L_t)_{t \in\PP^1}$ of
hyperplanes in~$\PP^{n+1}$ with  the base locus $M$ transversal
to $V$. Denote by $t_1$, \dots,~$t_N$ the collection of those~$t$
for which $L_t\cap V$ has singularities. Let $t_0$ be different
from either of $t_1$, \dots, $t_N$. Let $\gamma_i$~be a good
collection of paths in $\PP^1$
based in $t_0$. Let $e\in M-M \cap V$ be a base point. Let
${\cal V}_i$ be the variation operator 
corresponding to $\gamma_i$. Then the inclusion induces an
isomorphism:
\begin{equation}
 \pi_n(\PP^{n+1}-V,e) \longleftarrow \pi_n(L_{t_0}-L_{t_0}\cap V,e)\Big/
 \sum_{i=1}^N {\cal V}_i
\end{equation}
\end{theo}

There is affine version of this theorem equivalent to this 
one since $\pi_n(\PP^{n+1}-V)=\pi_n(\CC^{n+1}-V)$ in the case
when $H$ is transversal to $V$.

Recently, Cheniot and Eyral proposed definition of homotopy 
variation operator in general showed that the map
as in the above theorem is surjective (cf. \cite{chenioteyral};
see also \cite{tibar} for another discussion of variation operators).

\section{Local multivariable Alexander invariants: topological theory}

Now we want to develop an abelian version of the cyclic
theory presented so far. Though our goal at this point, as in the link theory,
is to study abelian covers, what will follow deviates from the 
link-theoretical point of view at several 
points. The most important one is that the Alexander type 
invariants are not polynomials. 
The substitute for the orders of the modules over PID which were 
discussed before are the subvarieties of commutative algebraic groups
called the characteristic varieties.

\subsection{Characteristic varieties of groups.}

\subsubsection{Definitions.}

Let us start with a classical construction of commutative algebra.
Let $R$ be a 
Noetherian commutative ring with a unit 
and let $M$ be a finitely generated 
$R$-module. Let  the homomorphism $\Phi: R^m \rightarrow R^n$ be such that
$M=\Coker \Phi$.  The  $k$-th Fitting ideal of $M$ is 
the ideal ${\cal F}_k(M)$ 
generated by $(n-k+1) \times (n-k+1)$ minors of the matrix 
of $\Phi$. ${\cal F}_k(M)$ 
depends only on $M$ rather than on $\Phi$.
The $k$-th characteristic variety $M$ is the reduced sub-scheme of 
$\Spec R$ defined by ${\cal F}_k(M)$. 
\par If $R={\bf C}[H]$ where $H$ is an abelian group then $\Spec R$ 
is a torus having the dimension equal to the rank of $H$.
If $H$ is free then after a choice of generators of $H$, $R$ can be 
identified with the ring of Laurent polynomials and 
$\Spec R=({{\bf C}^*})^{{\rm {rk}} H}$ is a complex
torus. In particular each $k$-th characteristic variety of an $R$-module 
is a subvariety 
$V_k(M)$ of $({{\bf C}^*})^{{\rm {rk}} H}$. If $H$ has a torsion 
then the number of connection components of $\Spec \CC[H]$ is the order of the 
torsion and the connected component of the unit can be identified with
$\Spec\CC[H/Tor(H)]$.

A more functorial description is the following (cf. \cite{BE}):
\begin{equation}
V_k(M)=\Supp_{red} (\wedge^k M)=\Supp_{red}(R/F_k(M))
\end{equation}

We shall apply this construction to the modules 
$A(G,\phi)$ defined in section \ref{alexsinvariants}. 
for pairs $(G,\phi)$ where $\phi: G \rightarrow \ZZ^r$.
Prime examples  which we shall consider are the following:

\begin{exam}
\par \noindent (i) Links in $S^3$.  
In this case $H_1(S^3-L,\ZZ)=\ZZ^r$ where $L$ is such a link
and $r$ is the number of its components.

\par \noindent (ii) Algebraic curves in $\CC^2$ having 
$r$ irreducible components (cf. section \ref{homologysection})
\end{exam}

We shall denote the corresponding 
characteristic varieties as $V_k(G,\phi)$ omitting 
$\phi$ when no confusion is possible.

\begin{dfn}(cf. \cite{charvar})
The {\it depth}
of a component $V$ of a characteristic variety $V_k(G)$
is the integer 
$i={\rm max} \{j \vert V \subset V_j(G)\}$.
\end{dfn}

In the case $r=1$ and $G$ is one of the groups as above, 
$V_1(G)$ is the zero set of the Alexander polynomial
and $V_k(G)$ is determined by the zero sets of elementary divisors of the 
Alexander module.
Vice versa, the zero sets of Fitting ideals determine
the zero sets of elementary divisors of a module over PID.
Since the orders of 
$\QQ[t,t^{-1}]$-modules in a cyclic decomposition 
 are determined up
to a unit of the ring of Laurent polynomials 
the depth of each root of the Alexander polynomial, 
given in terms of $V_k$'s, determines the Alexander module 
completely.

If ${\rm codim}V_1(G,\phi)$ in $\Spec \CC[G/G']$ is equal to 
one then the information carried by 
$V_1$ is equivalent to the multivariable Alexander 
polynomial up to the exponent of each factor 
(this is the case when $\phi$ is the abelianization 
of a link group). On the other hand if the codimension is bigger than 
one then for the 
pair $(G,\phi)$ the Alexander polynomial is not defined (or is trivial
depending on convention) but $V_1(G)$ can be very interesting.

Now, as the first example,
 let us calculate the characteristic varieties of a free group.
If $G=F_r$ is a free group on $r$-generators then $G'/G''
=H_1(\widetilde {\bigvee_r S^1},{\bf Z})$, where $\widetilde {\bigvee_r S^1}$
is the universal abelian cover of the wedge of $r$ circles.  
It fits into the exact sequence: $$0 \rightarrow 
H_1(\widetilde {\bigvee_r S^1},
{\bf C}) \rightarrow {\bf C}[{\bf Z}^r]^r \rightarrow I \rightarrow 0$$  
with $I$ denoting the augmentation ideal of the group ring of ${\bf Z}^r$:
$I=\Ker {\bf C}[\ZZ^r] \rightarrow \CC$
where the homomorphism sends each generator to $1 \in \CC$.
Indeed, as an universal abelian cover of 
$\widetilde {\bigvee_r S^1}$ one can take
the subset of ${\bf R}^r$ of points having at least $r-1$ integer
coordinates with the action of ${\bf Z}^r$ given by translations; unit 
vectors of the standard basis provide identification of 1-chains 
with ${\bf C}[{\bf Z}^r]^r$ while the module of 0-chains is identified 
with ${\bf C}[{\bf Z}^r]$. The boundary map sends each generator 
$e_i, i=1,...,r$ of 
${\bf C}[{\bf Z}^r]^r$ to $(t_i-1) \in {\bf C}[{\bf Z}^r]$.
This is the map which also appears in the Koszul complex 
(cf. \cite{localgebra}) in which we put $R={\bf C}[{\bf Z}]$:
\begin{equation}\label{koszul}
    \wedge^i R^r \longrightarrow \wedge^{i-1} R^r \longrightarrow
.... \rightarrow R^r 
\rightarrow R 
\end{equation}
where $\partial_i(e_{j_1} \wedge ... \wedge e_{j_i})=\sum (-1)^k(t_k-1)
e_{j_1} \wedge ... \wedge \hat e_{j_k} \wedge ... \wedge e_{j_i}$. 
Since $(t_1-1,...,t_r-1)$ is a system of parameters the complex 
(\ref{koszul}) is exact. Therefore
\begin{equation}
H_1(\widetilde{\bigvee_r S^1},{\bf C})=\Coker 
\Lambda^{r \choose 3}{\bf C}[{\bf Z}^r]^r \rightarrow 
\Lambda^{r \choose 2} {\bf C}[{\bf Z}^r]^r
\end{equation}
 in the Koszul resolution
corresponding to the $(t_1-1),...,(t_r-1)$.  
This implies that 
$V_i(F_r)={{\bf C}^*}^r$ for $0 < i \le r-1$ and $V_i(F_r)=(1,...,1)$
for $r \le i \le {r \choose 2}$ i.e. ${\CC^*}^r$ is component having 
depth $r-1$, and $1\in {\CC^*}^r$ has depth $r \choose 2$.

For arbitrary group $G$, as was pointed out in earlier sections, 
the Fox calculus provides presentation 
for the extension of the homology of universal abelian cover by 
the augmentation ideal of the group ring of the covering group.
This is sufficient to determine the characteristic varieties outside 
of the identity character.

\subsubsection{Unbranched covering}

The homology of a cyclic unbranched covering $X_n$ of a CW-complex $X$
with $\pi_1(X)=G$ corresponding to the homomorphism $G \buildrel \phi 
\over \rightarrow \ZZ 
\rightarrow \ZZ/n\ZZ$ 
can be found
using Milnor's exact sequence (cf. \cite{infmilnor}) i.e. the homology sequence
corresponding to the exact sequence of chain complexes 
\begin{equation}
0 \rightarrow C_*(\tilde X) \buildrel {t^n-1} \over \longrightarrow 
 C_*(\tilde X) \rightarrow C_*(X_n) \rightarrow 0
\end{equation}
The induced homology sequence:
\begin{equation}\label{milnorseq}
\rightarrow H_1(\tilde X, \CC) \buildrel {t^n-1} \over \longrightarrow 
H_1(\tilde X,\CC)
\rightarrow H_1(X_n,\CC) \rightarrow \CC \rightarrow 0
\end{equation}
shows that ${\rm rk} H_1(X_n,\CC)={\rm rk Coker}(t^n-1) 
\vert_{H_1(\tilde X,\CC)}+1$.
In abelian case, to find the homology of the covering $X_{n_1,...,n_r}$
corresponding to the homomorphism $G \buildrel \phi 
\over \rightarrow \ZZ^r \rightarrow \oplus_{i=1}^{i=r} \ZZ/n_i\ZZ$ 
the Milnor's sequence  (\ref{milnorseq}) should be replaced by the 
five term exact sequence corresponding to the spectral sequence 
of the covering group $H=\Ker \ZZ^r \rightarrow \oplus_i \ZZ/n_i\ZZ$ 
acting on the covering space $\tilde X$
corresponding to the homomorphism $\phi$:
\begin{equation}
      H_p(\ZZ^r,H_q(\tilde X,\CC)) \Rightarrow H_{p+q}(X_{n_1,...,n_r},\CC)
\end{equation}
This exact sequence is
\begin{equation}
H_2(X_{n_1,...,n_r},\CC) \rightarrow 
H_2(H,\CC) \rightarrow 
H_1(\tilde X)_H \rightarrow H_1(X_{n_1,...,n_r},\CC) \rightarrow H_1(H,\CC) \rightarrow 0
\end{equation}
where for a $H$-module $M$, $M_H=M/I(H)M$ is the module of covariants 
($I(H)$, as above, is the augmentation ideal of the group ring of $H$).
This yields the following formula for the first Betti number
of abelian covers:

\begin{prop} Let $X_{n_1,...,n_r}$ be the finite unbranched abelian cover 
of a CW-complex $X$ as above which is the quotient of the infinite
abelian cover corresponding to $\phi: G \rightarrow \ZZ^r$. 
Let $V_i(G,\phi)$ be the characteristic varieties of $(G,\phi)$.
For $P \in {{\bf C}^*}^r$ let 
$f(P,{G,\phi})=\{max \quad i \vert P \in V_i
(G,\phi) \}$. Then 
$${\rm rk}H_1({X_{n_1,...,n_r}})=
r+\Sigma_{\omega_i^{n_i}=1, (\omega_{n_1},...,\omega_{n_r}) \ne
(1,...,1)} f((\omega_{n_1},..,\omega_{n_r}),(G,\phi)) $$ 
\end{prop}

\subsubsection{Homology of local systems}

Homology of rank one local systems also can be described in 
terms of the characteristic varieties. Such a local system is 
a homomorphism $\chi: G \rightarrow \CC^*$ i.e. a character of the 
fundamental group. There is a natural identification: 
$\Spec \CC[G/G']$ and $\Char G$. Moreover,  
$\Spec \CC[G/\Ker \phi]$ can be identified with
 the subgroup of $\Char G$ of characters
which can be factored through $\phi$. We shall denote as
$\tilde X_{G/G'}$ is the 
infinite cover corresponding to the subgroup $G'$.
The homology $H_1(X,\chi)$ of the local system $\chi$ 
where $\chi \in \Char \pi_1(X)$ 
is defined as the homology of the chain complex:
\begin{equation}
  \rightarrow C_i(\tilde X_{G/G'}) \otimes_{\CC[G/G']} \CC_{\chi} 
\rightarrow C_i(\tilde X_{G/G'}) \otimes_{\CC[G/G']} \CC_{\chi}
\end{equation}
where the chains $C_i(\tilde X_{G/G'})$ of the universal 
abelian covers are given the structure of $\CC[G/G']$-module
and $\CC_{\chi}$ is $\CC$ with the module structure given by the 
character $\chi$.

One has the following:

\begin{prop}(cf. \cite{eko},\cite{charvar}) 
If $\chi \ne 1$ then 
 $$H_1(X,\chi)=
H_1(\tilde X_{G/G'},{\bf C}) \otimes_{{\bf C}[H_1(X,{\bf Z})]} {\bf C}_{\chi}$$
In particular, $\chi \in \Char G, \chi \ne 1$ belongs to $V_k(G)$ if and only
if $H_1(X,\chi) \ge k$.  
\end{prop}

\subsection{Links of plane curves and multivariable Alexander polynomial}

For a link in $S^3$ with $r$ components the characteristic varieties 
are just affine subvarieties of the torus. An interesting problem
is the following:

\begin{problem} Which sequence of subvarieties can occur as $V_i(G)$
where $G=\pi_1(S^3-L)$ for some link in $S^3$.
\end{problem}

For the multivariable Alexander polynomial one has:

\begin{equation}
\Delta(t_1^{-1},...,t_r^{-1})=\Delta(t_1,...,t_r)
\end{equation}
(up to a unit of the ring $\ZZ[\ZZ^r]$ i.e. a factor 
$\pm  t_1^{a_1}...t^{a_r}_r$ where $a_i \in \ZZ$)

The characteristic varieties of links of algebraic singularities
are very special. Let us call a translated subgroup 
of $\CC^r$ a coset of a subgroup $\CC^s$. Such a ``subgroup'' 
is said translated by an element of a finite order if this coset 
has finite order in $\CC^r/\CC^s$. Using the fact that the link
of algebraic singularities are iterated torus links one can prove 
the following:
\begin{prop}(cf. \cite{alexhodge})
The characteristic varieties of 
algebraic links are unions of translated subgroups.
\end{prop}

For example the link of singularity $x^r-y^r=0$ has the 
Alexander polynomial $t_1 \cdot ... \cdot t_r=1$. The Alexander 
polynomial of $(x^2-y^3)(x^3-y^2)=0$ is $(t_1^2t_2^3-1)(t_1^2t_2^3-1)=0$
(cf. \cite{alexhodge}).

\subsection{Links of isolated non normal crossings}

Disjoint non intersecting spheres of dimension greater than one and
having codimension 2 in an 
ambient sphere never can form a link of an algebraic singularity.
There is nevertheless a {\it local} abelian analog of the local cyclic theory 
of the links of algebraic singularities. It appears when one 
looks at isolated non normal crossings (cf. \cite{innc}, \cite{alexme}) 

\begin{dfn}(cf. \cite{innc})
 Let $D_1,..,D_k$ be divisors of a complex manifold $X$ and $P \in 
D_1 \cap ... \cap D_k$. These divisors have a normal crossing at 
$P$ if there exist in a  neighborhood   of $P$ in $X$ and 
a system of complex analytic local coordinates $(z_1,...,z_{{\rm dim} X})$
in it
such that $D_i$ is in this neighbourhood given by the equations $z_i=0$. 
$D_1,...,D_k$ have an 
isolated non normal crossing at $P$ is there exist a ball $B_{\epsilon}$ 
in $X$ centered at $P$ having sufficiently small radius $\epsilon$
such that for any $Q \ne P$ in $B_{\epsilon}$ the divisors $D_i$ containing 
$Q$ form at $Q$ a divisor with normal crossings.
\end{dfn}
In particular each $D_i$ has at most isolated singularity at $P$. 
A more general case, when the ambient space $X$ is allowed to have 
 a singularity at $P$ is considered in \cite{alexme}. The theory we shall 
describe is invariant under analytic changes of variables so we can assume 
that $X=\CC^{n+1}$. The starting point is the following vanishing result:

\begin{theo}\label{main}(cf. \cite{innc})
 Let $X=\bigcup_{i=1}^{r} D_i \subset {\bf C}^{n+1}$ 
be a union of 
$r$ irreducible germs of hypersurfaces with normal 
crossings outside of the origin. If $n \ge 2$,  
then $$\pi_1((\partial B_{\epsilon}-
\partial B_{\epsilon} \cap X)={\bf Z}^r \ \ {\rm and} \ \
\pi_k((\partial B_{\epsilon}-
\partial B_{\epsilon} \cap X)=0 \ \ \ {\rm for} \ \ 2 \le k <n. $$ 
\end{theo}

In the case when $r=1$ this result follows from Milnor's 
fibration theorem and connectivity of Milnor fibers 
(cf. \cite{singularpoints}). In fact, the universal cyclic 
cover of the complement to a link of isolated hypersurface singularity $D$
is homotopy equivalent to the Milnor fiber $M_D$. In particular 
$\pi_n(\partial B_{\epsilon}-D \cap \partial B_{\epsilon})=
H_n(M_D,\ZZ)$. For general INNC the main invariant is 
$\pi_n(\partial B_{\epsilon}-\cup D_i \cap \partial B_{\epsilon})$.
This, as usual, is the module over 
$\ZZ[\pi_1(\partial B_{\epsilon}-\cup D_i \cap \partial B_{\epsilon}]
=\ZZ[t_1,t_1^{-1},...,t_r,t_r^{-1}]$. We shall call it the
{\it homotopy module of INNC}.
In the case $r=1$ this module 
structure is equivalent to the module structure on an abelian group
with an automorphism where the abelian group is the middle homology of 
the Milnor fiber and the automorphism is the monodromy operator.
Notice that in the case of normal crossing (i.e. when the singularity 
is absent), the universal abelian cover of 
$\partial B_{\epsilon}-\cup D_i \cap \partial B_{\epsilon}$ is contractible 
and all homotopy groups are trivial. 

\begin{dfn}(cf. \cite{innc}) $k$-th characteristic variety $V_k(X)$
of an isolated non-normal crossing
$X=\cup_{i=1}^{i=r} D_i$ is the subset 
in ${\rm Spec}{\bf C}[\pi_1(\partial B_{\epsilon}-\partial B_{\epsilon} \cap 
(\bigcup_{1 \le i \le r} D_i))]$
formed by the zeros 
of the $k$-th Fitting ideal of   
$\pi_n(\partial B_{\epsilon}-\cup D_i \cap \partial B_{\epsilon})$
\end{dfn}

Let us consider an example of a non normal crossing.
The simplest non trivial case is when the components are given by 
linear equations i.e are given in $\CC^{n+1}$ by the equation 
$l_1 \cdot ...\cdot l_r=0$, where $l_i$ are {\it generic} linear forms
(i.e. a cone over a generic arrangement of 
hyperplanes in $\PP^n$). Since the complement to a generic arrangement
of $r$ hyperplanes in $\PP^n$ has a homotopy type of $n$-skeleton of 
the product of $r-1$-copies of the circle $S^1$ (in minimal cell decomposition
in which one has $r-1 \choose i$ cells of dimension $i$) one can calculate 
the module structure on the $\pi_n$ of such skeleton as follows
(cf. \cite{innc}). The universal 
cover of this skeleton 
is obtained by removing the $\ZZ^{r-1}$ orbits of all open faces of 
a dimension greater than $n$ in the unit cube in $\RR^{r-1}$. Hence 
$\pi_n(\partial B_{\epsilon}-D)=
H_n(\widetilde{Sk_n((S^1)^{r-1})},\ZZ)$ 
($\widetilde{Sk_n((S^1)^{r-1})}$ is the universal cover of the $n$-skeleton).
The chain complex of the universal cover of $(S^1)^{r-1}$ can be 
identified with the Koszul complex of the group ring of 
$\ZZ^{r-1}=\ZZ^r/(1,...,1)$ 
(so that the generators of $\ZZ^{r}$ correspond to the
standard generators of $H_1(\partial B_{\epsilon}-D)$). The system of 
parameters of this Koszul complex is $(t_1-1,..,t_r-1)$.
Hence $H_n(\widetilde{\partial B_{\epsilon}-D},\ZZ)=
Ker \Lambda^nR \rightarrow \Lambda^{n-1}R$ 
where $R=\ZZ[t_1,..,t_r,t_1^{-1},...,t_r^{-1}]/(t_1 \cdot ...\cdot t_r-1)$. 
As a result, one has the following presentation:
\begin{equation}
\Lambda^{n+1}({\bf Z}[t_1,t_1^{-1},...,t_r,t_r^{-1}]/(t_1...,t_r-1)^r)
\rightarrow \Lambda^n({\bf Z}[t_1,t_1^{-1},...,t_r,t_r^{-1}]/(t_1...,t_r-1)^r)
\rightarrow 
\end{equation}
$$\pi_n({\bf C}^{n+1}-\bigcup D_i) \rightarrow 0 $$
In particular, the support of the $\pi_n$ is the subgroup 
$t_1 \cdot ... \cdot t_r=1$.

The relation between the characteristic varieties, the unbranched covering 
spaces and the local systems described in the case of links in $S^3$ 
extends to this situation as well. We have the following:

\begin{prop} \label{unbranched}(cf. \cite{innc})
(a) For each 
$P \in {\rm Spec}{\bf C}[\pi_1(\partial B_{\epsilon}-
\partial B_{\epsilon} \cap X)]$ let 
$$f(P,X)=\{ {\rm max} \  k  \ \vert P \in V_k(X) \}$$
Let $U_{m_1,...,m_r}$ be unbranched cover 
of 
$\partial B_{\epsilon}-\partial B_{\epsilon} \cap 
(\bigcup_{1 \le i \le r} D_i)$ corresponding 
to the homomorphism 
$\pi_1(\partial B_{\epsilon}-\partial B_{\epsilon} \cap 
(\bigcup_{1 \le i \le r} D_i))=
{\bf Z}^r \rightarrow $ 
$\oplus_{1 \le i \le r} {\bf Z}/m_i{\bf Z}.$ 
Then $$ {\rm rk}H_p(U_{m_1,...,m_r},{\bf C})=\Lambda^p({\bf Z}^r)
\ \ {\rm for} \ \ p \le n-1,$$
$${\rm rk}H_n(U_{m_1,...,m_r},{\bf C})
=\sum_{(...,\omega_j,...),\omega_j^{m_j}=1} f((...,\omega_j,...),
\bigcup_{1 \le i \le r} D_i)$$

\bigskip

(b) If $1 \ne \chi \in \Char 
\pi_1(\partial B_{\epsilon}-\partial B_{\epsilon} \cap 
(\bigcup_{1 \le i \le r} D_i)=
\Spec \CC[\partial B_{\epsilon}-\partial B_{\epsilon} \cap 
(\bigcup_{1 \le i \le r} D_i)]$ is a character of the fundamental group
then 
$$H_i(\partial B_{\epsilon}-\partial B_{\epsilon} \cap 
(\bigcup_{1 \le i \le r} D_i),\chi)=0 \ \ 1 \le i \le n-1$$
 $$H_n(\partial B_{\epsilon}-\partial B_{\epsilon} \cap 
(\bigcup_{1 \le i \le r} D_i),\chi)=
\pi_n(\partial B_{\epsilon}-\partial B_{\epsilon} \cap 
(\bigcup_{1 \le i \le r} D_i)) \otimes_{\ZZ} \CC
 \otimes_{{\bf C}[H_1(X,{\bf Z})]} {\bf C}_{\chi}$$

\end{prop}

Milnor theory \cite{singularpoints} is applicable to 
INNC as to any hypersurface and one can relate  
relate Milnor's invariants to the characteristic varieties discussed here.
We have the following:

\begin{prop}\label{comparison}(cf. \cite{innc})
The homology of the Milnor fiber $F_D$ of the INNC singularity $D$ 
is given by: 
$$ H_p(F_D,{\bf Z})=\Lambda^p ({\bf Z}^{r-1}) \ \ {\rm for} \  1 \le p < n$$
The action of the monodromy of this homology is trivial.
The multiplicity of $\omega \ne 1$ as a root of the 
characteristic polynomial $\Delta_n(D,t)$ 
 of the Milnor's monodromy on $H_n(F_D,{\bf C})$ is equal to: 
$$m_{\omega}=f((...,\omega,...),D)=
{\rm max} \{ i \vert (\omega, ....,\omega) \in V_i(D)
 \subset {\rm Spec} 
{\bf C}[\pi_1 (\partial B_{\epsilon}-\partial B_{\epsilon} \cap D)] \}$$
\end{prop}

In the case of INNC the unbranched covering admits a natural compactification 
which provides model for the abelian branched covering 
of the sphere $S^{2n+1}$ with the link of INNC as the branching locus.
The branching cover itself is a link of an isolated complete 
intersection singularity. 

If the local equations of the locally irreducible components $D_1,...,D_r$
are $f_1,...,f_r$, then we can use as a model of abelian branched cover the 
link of singularity: 

\begin{equation}\label{icis}
z_1^{m_1}=f_1(x_1,...,x_{n+1}),...., z_1^{m_1}=f_1(x_1,...,x_{n+1}) 
\end{equation}

As a link of ICIS the link of singularity (\ref{icis}) is a $(n-1)$-connected
manifold having the dimension equal to $2n+1$. We shall express the 
homology of this link in terms of the characteristic varieties of the homotopy
modules associated to INNCs formed by various components of 
$D=\bigcup D_i$.

\begin{prop}\label{branchedcovers} 
Let $V_{m_1,...,m_r}$ be the link of singularity (\ref{icis}) 
which is the branched cover 
of $\partial B_{\epsilon}$ branched over \newline 
$\partial B_{\epsilon} \cap 
(\bigcup_{1 \le i \le r} D_i)$ with the Galois group 
$G=\oplus_{1 \le i \le r} {\bf Z}/m_i{\bf Z}$. 
For each 
$\chi \in {\rm Char} G$ let $I_{\chi}=
\{i \vert 1 \le i \le r, \chi({\bf Z}_{m_i}) \ne 1 \}$ 
where ${\bf Z}_{m_i}$ is the $i$-th summand of $G$. 
Any $\chi$ can also be considered as a character of 
$\pi_1(\partial B_{\epsilon}-\partial B_{\epsilon} \cap 
(\bigcup_{i \in I_{\chi}} D_i))$ 
in which case it will be called reduced and denoted $\chi_{red}$.
Let $V_{\chi}$ be the branched cover of $\partial B_{\epsilon}$
branched over $\partial B_{\epsilon} \cap 
(\bigcup_{i \in I_{\chi}} D_i)$ and having ${\rm Im \chi}=G/{\rm Ker} \chi$
as its Galois group.
Then $$ \pi_p(V_{m_1,...,m_r})=0 \ {\rm for} \ \ 1 \le p \le n-1,$$
$${\rm rk}H_n(V_{m_1,...,m_r},{\bf C})
=\sum_{\chi \in {\rm Char} } f(\chi_{red},
\bigcup_{i \in I_{\chi}} D_i).$$
\end{prop}
This proposition shows that there is a close relation between 
the homology of the tower  
of abelian covers and the characteristic 
varieties (at least in local case).
We shall use it in the following section 
for the calculation of the homology of infinite abalian covers in 
terms algebro-geometric data such as resolution of singularities
and the ideals of quasiadjunction.

\section{Hodge decomposition of local Alexander invariants.}
\label{localalexander}

\subsection{Zeros of Fitting ideals and Hodge numbers 
in cyclic case.}\label{zerosoffitting}

Our goal in this section is to describe the
structure of the characteristic varieties in 
the local case of the Alexander invariants of the 
germs of plane curves and for germs of INNC in terms of 
resolutions of singularities. This will give
an algebro-geometric description of these invariants. 
The global counterparts of the local invairants 
from this section will be considered in the section \ref{globalalexander}.
All structures introduced in this section are essentail 
for describing the global case. 

First let us consider the relationship between the 
Hodge structure of the cohomology of Milnor fiber and 
the Alexander invairants of the link of an isolated singularity.
In the cyclic case, the calculation of the Alexander polynomial
does not require Mixed Hodge theory and 
is a special case of  A'Campo's 
formula for the zeta function of the monodromy of a resolution
(\cite{A'Campo}).
Indeed, if $D$ has only one component with isolated singularity, 
the order of $\pi_n(S^{2n+1}-D \cap S^{2n+1}) \otimes \QQ$ 
is equivalent to the zeta function of the monodromy.
Hence, if $E_i$ are the components of the exceptional set
of a resolution $\pi$, $\pi^*(D)$ has the multiplicity $m_i$
along a component $E_i$ of the exceptional locus 
and the euler characteristic of the set of points in $E_i$ 
non-singular in the union of $\bigcup E_i$ and the proper preimage of $D$ 
is $\chi(E_i^{\circ})$ then (cf. \cite{A'Campo}): 
\begin{equation}
\Delta_n((S^{2n+1}-D \cap S^{2n+1})^{(-1)^n}(1-t)=
\Pi (1-t^{m(E_i)})^{-\chi(E_i^{\circ})}
\end{equation}
However even in cyclic case calculation of the zeros of higher 
Fitting ideals requires the mixed Hodge theory. We refer to 
\cite{deligne} or \cite{dimca} for the formalism of this structure.

In the rest of this section \ref{zerosoffitting} we shall 
focus mainly on the case of curves i.e. ${\rm dim}D=1$.
The cohomology 
group $H^1$\footnote{$^*$ we shall work with the cohomology as is more
common in Hodge theory though one has the dual structures on homology.
One of the differences is the presence of the negative weights in MHS 
in homology. See 
\cite{morgan}  
where the author works with MHS on homotopy groups (also discussion below
of the homotopy groups)
having negative weights and 
where natural dual theory with positive weights is not available}
of the Milnor fiber $(^*)$ of a plane curve singularity
supports a mixed Hodge structure with weights 0,1 and 2, 
with the identification 
\begin{equation}\label{logarithm}
N: W_2/W_1 \rightarrow W_0
\end{equation} 
given by the logarithm of an appropriate power of 
the monodromy(cf. \cite{oslo}). Recall that this means that 
one has canonically $(^{**})$ \footnote{$^{**}$ i.e. just an 
algebraic class of the germ at the zero
of fibration $\CC^2 \rightarrow \CC$ given by $(x,y) \rightarrow f(x,y)$}
defined (weight) 
filtration $H^1=W_2 \supseteq W_1 \supseteq W_0 \supseteq 0$ such 
that each quotient $W_n/W_{n-1}=\oplus_{p+q=n}H^{p,q}$. In fact there is 
a strong
relation between these groups $H^{p,q}$: they all come from increasing 
Hodge filtration. Moreover, if $T$ is the monodromy operator on $H^1(M,\CC)$ 
and $T=T_sT_u$ is the factorization into semisimple and unipotent part 
and if $N=log(T_u)=\sum_{i \ge 1} (-1)^{i-1}{{(T_u-I)^i} \over i}$ then 
$N$ induces the isomorphism in (\ref{logarithm}). 

All Hodge groups are invariant under the action of the 
semisimple part $T_s$ of the monodromy. 
Let $h^{p,q}_{\zeta}$ (cf. \cite{oslo}) be the dimension of the eigenspace
of this semisimple part acting on the space $H^{p,q}$.
The numbers $h^{p,q}_{\zeta}$
determine the Jordan form of the monodromy as follows.
The size of the Jordan blocks 
of the monodromy does not exceed 2 
and the number of blocks corresponding to
an eigenvalue $\zeta$ of size $1 \times 1$ (resp. $2 \times 2$)
is equal to $h^{1,0}_{\zeta}+h^{0,1}_{\zeta}$ (resp. $h^{0,0}_{\zeta}$). 
As a consequence, the generators of the Fitting ideals have the 
form:  
$$\Delta_i=\prod_{(\zeta)}(t-\zeta)^{a_{\zeta,i}}$$
where 

$$a_{\zeta,i}=\cases { h_{\zeta}^{1,0}+h_{\zeta}^{0,1}+2h_{\zeta}^{0,0}
-2(i-1) 
&  if  $1 \le i \le h_{\zeta}^{0,0}$  \cr
h_{\zeta}^{1,0}+h_{\zeta}^{0,1}-(i-1-h_{\zeta}^{0,0}) 
&  if  $ h_{\zeta}^{0,0} < i \le h_{\zeta}^{0,0}+h_{\zeta}^{1,0}+
h_{\zeta}^{0,1}$  \cr
0 &  if $i >  h_{\zeta}^{0,0}+h_{\zeta}^{1,0}+
h_{\zeta}^{0,1}$  \cr }$$
In particular,  $\Delta_i$  can be calculated
algebraically in terms of a resolution of the singularity 
since all Hodge numbers $h^{p,q}_{\zeta}$ can be found 
in terms of a resolution (cf. \cite{oslo}). 
A calculation of the Hodge numbers $h^{p,q}_{\zeta}$
 is equivalent to the identifying the following 
subsets in the set of zeros 
of the characteristic polynomial of the monodromy operator:
\begin{equation}\label{myspectrum}
{\cal H}_{p,q,k}=\{\zeta \vert h^{p,q}_{\zeta} \ge k \}
\end{equation}
Arnold-Steenbrink spectrum \cite{oslo} is also equivalent to this data.
Our goal for the rest of this section \ref{localalexander} will be to 
describe the partition of (the unitary part of) the 
zero sets of Fitting ideal i.e. the 
characteristic varieties of plane curve singularities (and also INNC) 
into sets (\ref{myspectrum}). We call this partition {\it the Hodge 
decomposition of characteristic varieties}.

In the abelian case the multivariable Alexander polynomial and hence
$V_1(G)$ can be found from a resolution as well (cf. \cite{EN}).

\begin{theo} Let $f_i=0, i=1,...r$ be the equations 
of branches of a reducible curve $C$. Let $\pi: \tilde \CC^2 \rightarrow 
\CC^2$ be a resolution of singularities and $E_j, j=1,...,N$ be 
the exceptional components. Let $m_{i,j}={\rm ord}_{E_j} \pi^*(f_i)$
and let $E_j^{\circ}$ be a Zariski open subset in $E_j$ consisting 
of points which are non singular on the reduced total preimage of $C$. 
Then the Alexander polynomial of the link of 
singularity of $C$ is:
$$\Pi_j (1-t_1^{m_{1,j}} \cdot ... \cdot t_r^{m_{r,j}})^{-\chi(E_i^{\circ})}$$
\end{theo}
For example for the Hopf link with $r$-components we obtain 
$(1-t_1 \cdot...\cdot t_r)^{r-2}$.

\bigskip

We shall use the description of the cohomology of 
branched coverings from the last lecture to calculate higher $V_i(G)$. 
The information about 
the characteristic varieties is closely related to the 
information about the cohomology of branched abelian covers (cf. Prop.
\ref{branchedcovers}). Those 
are the links of singularities of complete intersection and 
hence have the mixed hodge structure. We shall calculated 
the eigenspaces of the deck transformations acting on the Hodge 
spaces of the cohomology of branched coverings 
in  terms of the algebraic data of the singularity i.e. some
associated ideals generalizing the ideals of quasiadjunction
defined earlier in this paper. This will 
give us the calculation of the sets (\ref{myspectrum}) and hence 
the characteristic varieties. The MHS on the link can be 
described in terms of resolution of the singularity so we shall need 
a description of the forms of resolution of singularities which we
shall describe in the next section. We shall start it by description 
of a global precursor: adjoint hypersurfaces. 

\subsection{Theory of adjoints}\label{adjointsection} 

The classical theory of adjoints (which already was briefly mentioned
in section \ref{alexpolynomialposition} and which was used there 
to calculated one variable Alexander polynomial)
gives a presentation of the  
geometric genus of a resolution of singularities of
a hypersurface or a complete intersections in $\PP^N$
in terms of the degree of hypersurface 
(or the multidegree in the case of a CI) and 
local data about the singularities. A starting point maybe 
an observation that the genus of a plane curve $C$ of degree $d$ (i.e. 
$h^{1,0}=\dim H^0(\Omega_C)$) is equal to ${{(d-1)(d-2)} \over 2}$ which 
also is equal to the 
dimension of plane curves of degree $d-3$. If $C$ is not smooth 
but have $\delta$ nodes then the genus of desingularization 
${{(d-1)(d-2)} \over 2} -\delta$ is 
the  the dimension of space of plane curves of degree $d-3$ passing through 
the nodes. In the case when $C$ has singularities more complicated than nodes 
one can associate with each singular point $P$ the ideal
in the local ring $\O_P$ (adjoint ideal) 
such that the genus of desingularization 
is the dimension of the space of curves of degree $d-3$ passing through 
the singularities of $C$ and which local equations at $P$ belong to the 
adjoint ideal in $\O_P$.

Explicitly, 
these ideals can be described as follows. Let $X \subset \PP^{n+1}$ be 
a hypersurface and 
let $f_*: \tilde X \rightarrow X$ 
be a resolution of singularities. 
Let ${\cal A}=f_*(\Omega^n_{\tilde X})(-d+n+2)$
and ${\cal A}'=\pi^{-1}({\cal A})$ 
where $\pi: \O_{\PP^{n+1}} \rightarrow \O_X$ 
is the restriction map. Then $\cal A'$ is called the sheaf of adjoint ideals 
and for example $H^0({\cal A}'(d-n-2))=h^{n,0}(\tilde X)$.

The sheaf $f_*(\Omega^n_{\tilde X})$ is a subsheaf 
of $i_*(\Omega^n_{X-\Sing X})$ where $i: X-\Sing X \rightarrow X$ is the 
embedding. One has the residue map which fit into the exact 
sequence:
\begin{equation}
0 \rightarrow \Omega^{n+1}_{\PP^{n+1}} \rightarrow 
\Omega^{n+1}_{\PP^{n+1}}(X) \buildrel {\rm Res} \over 
 \rightarrow i_*(\Omega^n_{X-\Sing X}) \rightarrow 0
\end{equation}
The residue map sends a form 
$\omega={{f(z_1,...,z_{n})dz_1 \wedge ... \wedge dz_{n}} \over F}$, 
defined in a chart with coordinates $z_1,...,z_n$ and
having the pole of order one 
along $X$ (given by in this chart by the equation $F=0$)
to $(-1)^{j-1}{{fdz_1 \wedge ...\hat dz_j \wedge dz_{n}} 
\over {F_{z_j}}}\vert_{X-\Sing X}$ (the restriction is 
independent of $j, 1 \le j \le n$).
From this point of view the stalk of the sheaf $\cal A'$ at $P \in \PP^{n+1}$ 
consists of $f \in \O_P$ such that $Res (\omega)$ extends to a holomorphic form
on some resolution $\tilde X$. 

If the case when a subvariety $X \subset \PP^{n+r}$ 
is a complete intersection of 
hypersurfaces $F_1=...=F_r$ the stalk of the sheaf of adjoint ideals at a 
singular point can be described 
using the fact that a holomorphic $n$-form can be obtained as a multiple 
of the residue of $(n+r)$-form
having poles of order one along each hypersurface $F_i=0$ and that 
\begin{equation}\label{icisresidue}
Res {{dw_1 \wedge .... \wedge dw_{n+r}} \over {F_1 \cdot ....\cdot F_r}}=
{{... \wedge \hat dw_{i_1}... \wedge \hat dw_{i_r} ...} \over
 {{\partial (F_1,...,F_{r})} \over
 {\partial (w_{i_1},...,w_{i_r})}}}\vert_X
\end{equation}
where ${{\partial (F_1,...,F_{r})} \over
 {\partial (w_{i_1},...,w_{i_r})}}$ is the Jacobian of partial derivatives 
of the system $(F_1,...,F_r)$ relative to variables $(w_{i_1},...,w_{i_r})$
(it is easy to check that the restriction up to sign is independent of 
collection of variables $(w_{i_1},...,w_{i_r})$).

This construction describes the differential forms 
on a resolution in terms of the linear systems 
of hypersurfaces in $\PP^{n+r}$ given by ideal sheaves on the latter. 
The cohomology and the Hodge structure on the link of a complete 
intersection singularity can be described in terms of differential 
forms (cf \cite{oslo}) and this description 
can be used to calculate the cohomology of the link 
in terms of certain ideals in the local ring of singular point.
We shall review the construction of MHS on various cohomology associated with 
singular points of complex spaces and in the next section we shall 
discuss the connections
with the ideals of quasiadjunction.

The cohomology of the link $L$ of an isolated singularity $x$ of a 
complex space $X$ (${\rm dim}X=n$) can be given a Mixed Hodge 
structure, for example 
using canonical identification $H^k(L)=
H^*_{\{x\} }(X)$ with the local cohomology ($^*$)
\footnote{* Recall that if $Y$ is a subset in a topological space $X$ 
and $\cal F$ is a sheaf on $X$ 
then $H^i_Y(X)$ is the right derived functor of the functor 
$\Gamma_Y(X,{\cal F})$ 
of sections of $\cal F$ supported on $Y$. It fits into long 
exact sequence: 
 $ ... \rightarrow H_Y^i(X,{\cal F}) \rightarrow H^i(X,{\cal F})
\rightarrow 
H^i(X-Y,{\cal F}) \rightarrow  H_Y^{i+1}(X,{\cal F}) \rightarrow ...$
(cf. \cite{hartshorne})}.
The mixed Hodge structure on the latter was described 
in \cite{arcata}. 
The Hodge numbers: $h^{kpq}(L)=
{\rm dim} Gr^p_FGr^W_{p+q}H^k(L)$ have the following 
symmetry properties:
\begin{equation}
h^{kpq}=h^{2n-k-1,n-p,n-q}
\label{symmetry}
\end{equation} 
If $E$ is the exceptional divisor for a resolution, then
for $k <n$ one has  
$$
h^{kpq}(L)=h^{kpq}(E) \ \ \ if \ p+q <k$$
\begin{equation}
h^{kpq}(L)=h^{kpq}(E) - h^{2n-k,n-p,n-q}(E)
\ \ \ if \ p+q =k
\label{reduction1}
\end{equation} 
$$h^{kpq}(L)=0 \ \ \ if \ p+q >k$$
\par The local 
cohomology $H^*_E(\tilde X)$ (\cite{arcata}) where $\tilde X$ is a 
resolution of $X$ support the canonical mixed Hodge structure. 
Let us consider it in more detail 
in the case $dim_{\bf C} \tilde X=2$ which we shall need to 
describe the characteristic varieties in the case of germs of plane curves. 
We have 
\begin{equation}
H^*_E(\tilde X)=Hom(H^{4-*}(E),{\bf Q}(-2))
\label{duality}
\end{equation}
where ${\bf Q}(-2)$ is Tate Hodge of type $(2,2)$. 
Since the Hodge and weight filtrations on 
$H^1(E)$ have the form: 
$$H^1(E)= W_1 \supset W_0 \supset 0, \ \ \ H^1(E)=F^0 \supset F^1 \supset
F^2=0$$
on $H^3_E(\tilde X)$ we have:
 
$$H^3_E(\tilde X)=W_4 \supset W_3 \supset W_2=0,H^3_E(\tilde X)=F^1 \supset F^2 \supset
F^3=0$$

Moreover 
\begin{equation}
F^1H^1(L)=F^1H^1(E)=F^2H^3_E(\tilde X)
\label{hodgefiltration}
\end{equation}

In order to describe this 
mixed Hodge structure one can use the following complex: 

\begin{equation}
0 \rightarrow A^2_E(\tilde X) \rightarrow A^3_E(\tilde X) 
 \rightarrow 0
\label{complex}
\end{equation}

where $$A^2_E(\tilde X)=
\Omega^1_{\tilde X}(log \ E)/\Omega^1_{\tilde X},
A^3_E(\tilde X)=
\Omega^2_{\tilde X}(log \ E)/\Omega^2_{\tilde X}$$

with filtrations given by 
$$F^2A^p_E(\tilde X)=0 \ for \ p<3,
F^2A^p_E(\tilde X)=A^p_E( \tilde X) \ for \ p \ge 3$$

$$W_3A^3_E(\tilde X)=W_1\Omega^2_{\tilde X}(log \ E)/\Omega^2_{\tilde X}$$

Since $H^3(E)=0$, the relations (\ref{reduction1}) and (\ref{duality}) yield 
that the complex (\ref{complex}) completely determines 
$h^{1pq}$ (and hence all Hodge numbers $h^{kpq}$ by (\ref{symmetry})).

Putting all this together we obtain the following isomorphism:

\begin{equation}\label{keyexpression}
H^0(\Omega_{\tilde X}^2(log \ E))
/H^0(\Omega_{\tilde X}^2)=F^1 H^1(L)
\end{equation}

Our next goal will be to apply this to the link of complete intersection 
singularity which is the abelian cover branched over the 
a link plane curve singularity i.e. to the link (\ref{icis})
which in case of curves with $r$ components is the following link in 
$\CC^{r+2}$:

\begin{equation}\label{iciscurves}
 z_1^{m_1}=f_1(x,y),\ \ ... \  \  z_1^{m_1}=f_r(x,y)
\end{equation}

We want to calculate the eigenspaces corresponding to the characters
of the Galois group acting on the Hodge spaces $H^{p,q,k}$ of 
singularity (\ref{iciscurves}) by interpreting the left hand side of 
(\ref{keyexpression}) in terms of ideals in the local ring of the 
singularity $f_1\cdot \cdot \cdot f_r=0$ in $\CC^2$.

\subsection{Ideal of quasiadjunction and log-quasiadjunction}

We shall start with the following multivariable generalization
of the ideals of quasiadjunction introduced in section
\ref{alexpolynomialposition}.
\begin{dfn} \label{formquasiadjunction}
(cf. \cite{wall} \cite{charvar})
{\it An ideal of quasiadjunction} 
of type $(j_1,..,j_r \vert m_1,..,m_r)$ 
is the ideal in the local ring of the singularity of
$C$ (i.e. $O \in {\bf C}^2$) consisting of germs $\phi$ such that 
the 2-form:
$$
\omega_{\phi}={{\phi(x,y) z_1^{j_1} \cdot \cdot \cdot z_r^{j_r} dx  \wedge 
dy} \over  {z_1^{m_1-1} \cdot \cdot \cdot z_r^{m_r-1}}},
\label{form}$$
extends to a holomorphic form on a resolution of the singularity 
of the abelian cover of a ball $B$ of type $(m_1, ...., m_r)$,
i.e. a resolution of (\ref{iciscurves})
(we suppress dependence of $\omega_{\phi}$ on $j_1,..j_r, m_1,...m_r$).
In other words, $\phi z_1^{j_1} \cdot \cdot \cdot z_r^{j_r}$ belongs to 
the adjoint ideal of the singularity (\ref{iciscurves}).
In particular the condition on $\phi$ is independent of resolution.
\end{dfn}
Note that $\omega_{\phi}$ in definition (\ref{formquasiadjunction}) 
 is the residue of the form 
${{\phi z_1^{j_1} \cdot \cdot \cdot z_r^{j_r} dx  \wedge 
dy} \over  {(z_1^{m_1}-f_1(x,y)) \cdot \cdot \cdot (z_r^{m_r}-f_r(x,y))}}$
(cf. (\ref{icisresidue}))
We always shall assume that 
$0 \le j_1 < m_1,..,0 \le j_r < m_r$. Also, notice that forms 
$\omega_{\phi}$ in definition (\ref{formquasiadjunction}) 
are exactly the forms on the abelian cover 
(\ref{iciscurves}) which are the eigenforms corresponding to the character
of the Galois group taking value $exp({{2 \pi i(j_i-m_i+1)} \over {m_i}})$
 on the automorphism of the surface (\ref{iciscurves}) 
induced by multiplication of the $i$-th coordinate 
by $exp({{2 \pi i} \over {m_i}})$.  

{\it An ideal of log-quasiadjunction} 
(resp. {\it an ideal of weight one} {\it log-quasiadjunction}) 
of type $(j_1,..,j_r \vert m_1,..,m_r)$
is the ideal in the same local ring 
consisting of germs $\phi$ such that $\omega_{\phi}$ 
extends to a log-form (resp. weight one log-form) on a resolution 
of the singularity of the same abelian cover. Recall (cf. \cite{deligne})
that a holomorphic 2-form is weight one log-form if it is a 
combination of forms having
poles of order at most 
one on each component of the exceptional divisor and  
not having poles of order one on a pair of intersecting components. 
These ideals are also independent of a resolution (cf. \cite{alexhodge}). 

One can show (cf. \cite{charvar}) that an ideal of quasiadjunction 
${\cal A}(j_1,..,j_r \vert m_1,...,m_r)$ is determined by the vector
(i.e. depends only on the collection of ratios):
\begin{equation}
({{j_1+1} \over {m_1}},....,{{j_r+1} \over
  {m_r}}).
\label{vector}
\end{equation}
 This is also the case for the ideals of 
log-quasiadjunction and weight one log-quasiadjunction. 
Indeed, these ideals can be described in terms of resolutions 
as follows. 
For a given embedded resolution $\pi: V \rightarrow {\bf C}^2$ of the 
germ $f_1 \cdot \cdot \cdot f_r=0$ with the exceptional curves 
$E_1,..,E_k,...,E_s$ let $a_{k,i}$ (resp. $c_k$, resp. $e_k(\phi)$)
 be the multiplicity of the pull back on $V$ of 
$f_i$ ($i=1,..,r$) (resp. $dx \wedge dy$, resp. $\phi$)
along $E_k$.
Then $\phi$ belongs to the ideal of quasiadjunction of type 
   $(j_1,..,j_r \vert m_1,..,m_r)$ if and only if for any $k$ 
\begin{equation}
 a_{k,1}{{j_1+1} \over {m_1}}+....+a_{k,r}{{j_r+1} \over {m_r}}
 > a_{k,1}+...+a_{k,r}-e_k(\phi)-c_k-1
\label{quasiadj}
\end{equation}
(cf. \cite{charvar}).
Similar calculation shows that 
a germ $\phi$ belongs to the ideal of log-quasiadjunction
corresponding to 
$(j_1,..,j_r \vert m_1,..,m_r)$ if and only if the 
inequality 
\begin{equation}
a_{k,1}{{j_1+1} \over {m_1}}+....+a_{k,r}{{j_r+1} \over {m_r}}
 \ge  a_{k,1}+...+a_{k,r}-e_k(\phi)-c_k-1
\label{logquasiadj}
\end{equation}
is satisfied for any $k$.
In addition,  a germ $\phi$ belongs 
to the ideal of weight one log-quasiadjunction if and
only if this germ is a linear combination of germs $\phi$ 
satisfying inequality (\ref{logquasiadj})
for any collection 
of $k$'s such that corresponding components do not 
intersect 
and satisfying the inequality  (\ref{quasiadj}) 
for $k$ outside of this collection. 
We shall denote the ideal of quasiadjunction (resp. log-quasiadjunction, 
resp. weight one log-quasiadjunction) corresponding to 
$(j_1,..,j_r \vert m_1,..,m_r)$ as ${\cal A} (j_1,..,j_r \vert
m_1,..,m_r)$ (resp. ${\cal A}'' (j_1,..,j_r \vert
m_1,..,m_r)$, resp. ${\cal A}' (j_1,..,j_r \vert
m_1,..,m_r)$). Note the inclusions: 
$${\cal A} (j_1,..,j_r \vert
m_1,..,m_r) \subseteq  
{\cal A}' (j_1,..,j_r \vert
m_1,..,m_r) \subseteq 
{\cal A}'' (j_1,..,j_r \vert m_1,..,m_r)$$

\par \noindent
Both (\ref{quasiadj}) and (\ref{logquasiadj}) follow 
from the following calculation (cf. \cite{charvar} section 2
for complete details).
One can use the normalization of the fiber product
$\widetilde V_{m_1,...,m_r}=
V \times_{{\bf C}^2} V_{m_1,..,m_r}$
as a resolution of singularity (\ref{iciscurves})
in the category of manifolds with quotient singularities 
(cf. \cite{lipman}).
We have: 
\begin{equation}
\matrix{  \tilde V_{m_1,..,m_r} & \buildrel \tilde p \over 
                                                \rightarrow & V \cr
               \tilde \pi  \downarrow & &  \pi \downarrow \cr
               V_{m_1,..,m_r}&  \buildrel p \over 
                     \rightarrow & {\bf C}^2 \cr }
\label{resolution}
\end{equation}        
The preimage of the exceptional divisor of $V \rightarrow {\bf C}^2$
in $\widetilde V_{m_1,...,m_r}$ forms a divisor with normal crossings
(cf. \cite{oslo}), though the preimage of each component is reducible 
in general. In this case the irreducible components 
above each exceptional curve do not intersect.
If the Galois group $G$ of $\tilde p$ is abelian
(as we always assume here) and, in particular, is the quotient of 
$H_1(B-C \cap B,{\bf Z})$, 
then the Galois group of ${\tilde p}^{-1}(E_i) \rightarrow E_i$ 
is $G/(\gamma_i)$ 
where for an exceptional curve $E_k$, $\gamma_k$ is the image in 
the Galois group of the homology class of the boundary of a small 
disk transversal to $E_k$ in $V$. The components of 
${\tilde p}^{-1}(E_i)$ correspond to the elements of 
$G/(\gamma_i, ...\gamma_l ..)$ where $l$ runs through indices 
of all exceptional curves intersecting $E_i$, while 
$\tilde p_i$ restricted on each component has 
$(\gamma_i, ...\gamma_l ..)/(\gamma_i)$ as the Galois group.
The points ${\tilde p}^{-1}(E_i \cap E_j)$ 
correspond to the elements of $G/(\gamma_i,\gamma_j)$
and the points of ${\tilde p}^{-1}(E_i \cap E_j)$
belonging to a fixed component correspond to cosets
in $ (\gamma_i, ...\gamma_l ..)/(\gamma_i, \gamma_j).$
The order of the vanishing of $\omega_{\phi}$ on 
$\widetilde V_{m_1,...,m_r}$ along $E_k$ is equal to:
\begin{equation}
\Sigma_{i=1}^{i=r}
 ( j_i-m_i+1) {{m_1 \cdot \cdot \cdot \hat m_i \cdot \cdot \cdot m_r
 \cdot a_{k,i}} \over {g_{k,1} \cdot \cdot \cdot g_{k,r}s_k}}+
{{m_1 \cdot \cdot \cdot m_r \cdot ord_{E_k}(\pi^*(\phi))} \over {g_{k,1} \cdot
\cdot \cdot g_{k,r} \cdot s_k}}+
\label{orderofzero}
\end{equation}
$$+{{c_k \cdot m_1 \cdot \cdot \cdot m_r} \over {g_{k,1} \cdot \cdot \cdot
 g_{k,r} \cdot s_k}}+ {{m_1 \cdot \cdot \cdot m_r} \over {g_{k,1} \cdot
\cdot \cdot g_{k,r} \cdot s_k}}-1 $$
where $g_{k,i}=g.c.d.(m_i,a_{k,i})$ and 
$s_k=g.c.d.(...,{{m_i} \over {g_{k,i}}},...)$.

A consequence of (\ref{orderofzero}) is that $\omega_{\phi}$ has 
an order of pole equal to one (resp. zero) along the component $E_k$ of the 
above resolution if and only if for such $\phi$ one has 
equality in (\ref{logquasiadj}) (resp. (\ref{quasiadj}) is satisfied). 
\begin{prop}(cf. \cite{alexhodge})
 1. Let ${\cal A}''$ be an ideal of log-quasiadjunction. 
There is a unique polytope ${\cal P}({\cal A}'')$ such that
a vector $({{j_1+1} \over {m_1}},...,
 {{j_r+1} \over {m_r}}) \in {\cal P}({\cal A}'')$ 
if and only if  the ideal ${\cal A}'' (j_1,..,j_r \vert
m_1,..,m_r)$ contains ${\cal A}''$ 
\footnote{$(^*)${ i.e. a subset in ${\bf R}^r$ given by a set 
of linear inequalities $L_s \ge k_s$. 
We say that an affine hyperplane in ${\bf R}^r$ supports a codimension one face
of a polytope if the intersection of this hyperplane with the boundary 
of the polytope has dimension $r-1$. 
A face of a polytope is the intersection of a supporting face of 
the polytope with the boundary. A codimension one face of a polytope 
in ${\bf R}^r$ is a polytope of dimension $r-1$.
By induction one obtains faces of arbitrary codimension
for original polytope (for $r=3$ those are called edges and vertices). 
The boundary of the polytope is the union of its faces.}}.
\par 2. The set of vectors (\ref{vector}) for which 
${\cal A} (j_1,..,j_r \vert m_1,..,m_r) \ne$ 
\newline ${\cal A}''(j_1,..,j_r \vert m_1,..,m_r)$ is a dense 
subset in the boundary of the polytope having as its closure 
a union of faces of such a polytope.
The closure of the set of vectors (\ref{vector}) 
for which  ${\cal A}' (j_1,..,j_r \vert m_1,..,m_r) \ne 
{\cal A''}(j_1,..,j_r \vert m_1,..,m_r)$
is also a union of certain faces of such a polytope. 
\par 3. The ideal ${\cal A}(j_1,..,j_r \vert  m_1,...,m_r)$ 
 (resp. ${\cal A}'(j_1,..,j_r \vert  m_1,...,m_r)$ and 
\newline ${\cal A}''(j_1,..,j_r \vert  m_1,...,m_r)$) is independent of 
the array $(j_1,..,j_r \vert m_1,..,m_r)$ as long as the vector (\ref{vector}) 
varies within the interior of the same face of quasiadjunction.
\end{prop}

\noindent We shall call the above faces {\it the faces of quasiadjunction} 
(resp. {\it weight one faces of quasiadjunction}). ${\cal A}_{\Sigma}$ will 
denote ${\cal A}(j_1,...,j_r \vert m_1,...,m_r)$ with corresponding
vector (\ref{vector}) belonging to the interior of a face of quasiadjunction $\Sigma$
(similarly for ${\cal A}_{\Sigma}'$ and ${\cal A}_{\Sigma}''$).

In the case $r=1$ and when $f(x,y)$ is weighted homogeneous 
one can use the description of the adjoint ideals given by M.Merle 
and B.Tessier (cf. \cite{merle} and section \ref{alexpolynomialposition}). 
The polytopes of quasiadjunction 
are in $\RR$ and hence are just constants. They are 
the constants of quasiadjunction introduced in \cite{arcata1}. 
It was shown in \cite{loeservac} that they are the elements 
of Arnold-Steenbrink spectrum which belong to the interval (0,1).

The polytopes of quasiadjunctions are subsets of 
a unit cube $\cal U$ with the coordinates
corresponding to the components of the link.
We shall view it also 
as the fundamental domain for the Galois group 
$H^1(S^3-L,{\bf Z})$
of the universal 
abelian cover $H^1(S^3-L,{\bf R})$
 of the group $H^1(S^3-L,{\bf R}/{\bf Z})$
of the unitary characters of 
$H_1(S^3-L,{\bf Z})$
(i.e. the maximal compact subgroup 
of $Char (H_1(S^3-L,{\bf Z}))=H^1(S^3-L,{\bf C}^*)$).
$exp: {\cal U} \rightarrow 
Char (H_1(S^3-L,{\bf Z}))$ will denote the restriction of 
$H^1(S^3-L,{\bf R}) \rightarrow H^1(S^3-L,{\bf R}/{\bf Z})$
on $\cal U$.
\par For any sub-link $\tilde L$ of $L$, i.e. a link formed  by 
components of $L$, we have 
surjection $\pi_1(S^3-L) \rightarrow 
\pi_1(S^3-\tilde L)$ induced by inclusion. 
Hence $Char H_1(S^3-\tilde L,{\bf Z})$
is a sub-torus of $Char H_1(S^3-L,{\bf Z}))$ (in coordinates in the latter 
torus corresponding to the components of $L$ it is given by equations 
of the form $t_{\alpha}=1$ where subscripts correspond to components of $L$
absent in $\tilde L$). 
Moreover, since the homology of the universal abelian cover 
$H_1(\widetilde {S^3-L})$ surjects onto
$H_1(\widetilde {S^3-\tilde L})$,
 it follows that $V_i(S^3-\tilde L)$ belongs to a 
component of $V_i(S^3-L)$ (cf. \cite{charvar}). 
We shall call a character of $\pi_1(S^3-L)$ 
(or a connected component of $V_i(S^3-L)$) 
{\it essential} if it does not belong to a subtorus $Char H^1(S^3 -\tilde L)$ 
for any sublink $\tilde L$ of $L$.
\par Let $L_{m_1,..,m_r}$ be the link of singularity
(\ref{iciscurves}) or equivalently the cover 
of $S^3$ branched over the link $L$ and having a quotient 
$H_{m_1,..,m_r}={\bf Z}/m_1{\bf Z} \oplus .... \oplus 
{\bf Z}/m_1{\bf Z}$ of  $H_1(S^3-L,{\bf Z})$ 
as its Galois group. We shall view $Char H_{m_1,..,m_r}$
as a subgroup of $Char H_1(S^3-L,{\bf Z})$. The group 
$H_{m_1,..,m_r}$ acting on 
$H^1(L_{m_1,..,m_r})$ preserves both Hodge and weight filtrations.

\begin{theo}(cf. \cite{alexhodge}) An essential  character $\chi \in 
Char (H_1(S^3-L,{\bf Z})) $ is a character of 
the representation of $H_{m_1,..,m_r}$ acting on 
$F^1(H^1(L_{m_1,..,m_r}))$ if and only if it factors through
the Galois group $H_{m_1,..,m_r}$ and belongs to the image
of a face of quasiadjunction under the exponential map.
\par The multiplicity of $\chi$ in this representation of the
Galois group is equal to $dim {\cal A}_{\Sigma}''/{\cal A}_{\Sigma}$
where $  {\cal A}_{\Sigma}''$ (resp. ${\cal A}_{\Sigma}$)
is the ideal of log-quasiadjunction (resp. ideal of quasiadjunction)
corresponding to a vector (\ref{vector}) belonging to 
the face of quasiadjunction $\Sigma$. 
\par A character $\chi$ is a character of
the representation of the Galois group of the cover on
$W_0(H^1(L_{m_1,..,m_r}))$ if and only if it factors through
the Galois group $H_{m_1,..,m_r}$ and it belongs to the image
under the exponential map of a weight one face of quasiadjunction.
\end{theo}

\subsection{Multiplier ideals and log-canonical thresholds.}
\label{multidealsection}

The ideals and polytopes of quasiadjunction are closely related 
to recently studies multiplier ideals (cf. \cite{nadel}, \cite{lazarsfeld})
and log-canonical thresholds.

For a ${\bf  Q}$-divisor $D$ on a non singular
manifold $X$ its multiplier ideal  
${\cal J}(D)$ (cf. (\cite{nadel}) 
can be defined as follows. Let $f: Y \rightarrow X$ be an embedded 
resolution of $D$ and $f^*(D)=-E$. Then 
${\cal J}(D)=f_*({\cal O}_Y(K_Y-f^*(K_X)-\lfloor E \rfloor))$ 
where $\lfloor E \rfloor$
is round-down of a $\bf Q$-divisor. 
In this terminology one can define the ideals of 
quasiadjunction as follows. 
For an array $(\gamma_1,..,\gamma_r), (\gamma_i \in {\bf Q})$ 
let $D_{\gamma_1,..,\gamma_r}$ 
be given by equation $f_1^{\gamma_1} \cdot \cdot \cdot 
f_r^{\gamma_r}$. 
Then ${\cal J}(D_{\gamma_1,..,\gamma_r})={\cal A}(j_1,..,j_r \vert m_1,...,
m_r)$ where $\gamma_i=1-{{j_i+1} \over {m_i}}$ for $i=1,..,r$. 
This follows immediately from (\ref{quasiadj}).

To describe the relation with the log-canonical thresholds, 
recall (\cite{kollar})
that a pair $(X,D)$ where $X$ is normal and $D$ is 
a $\bf R$-divisor such that $K_X+D$ is $\bf R$-Cartier 
is called log-canonical at $x \in X$
if for any birational morphism
$f: Y \rightarrow X$, with $Y$ normal, in the decomposition 
\begin{equation} K_{Y}=f^*(K_X+D)+\sum_E a(E,X,D)E
\label{discrepancy}
\end{equation} 
for each irreducible $E$ having center at $x$ one has $a(E,X,D) \ge -1$. This 
coefficient is called {\it discrepancy} of divisor $D$ 
on $X$ along $E$.

\begin{prop}(cf. \cite{alexhodge}) The local ring ${\cal O}_O$  
of a singularity $f_1 \cdot \cdot \cdot f_r=0$ at the origin $O$
of ${\bf C}^2$ considered as the ideal in itself 
is an ideal of log-quasiadjunction.
Let $\cal P$ be the corresponding polytope of log-quasiadjunction.
Let $D_i$ be the divisor in ${\bf C}^2$ with the local equation $f_i=0$ 
near the origin. 
\newline \noindent Then for $\{(\gamma_1,...,\gamma_r ) \} \in  {\bf R}^r$
the divisor $\gamma_1 D_1+...+\gamma_r D_r$ is log-canonical 
at $(0,0) \in {\bf C}^2$ if and only if $(1-\gamma_1,..,1-\gamma_r)$
belongs to the polytope $\cal P$.
\label{treshold}
\end{prop}

To see why this is the case, let us consider the polytope given 
by inequalities (\ref{quasiadj}) in which one puts 
$e_k({\cal A}'')=0$, i.e. 
 \begin{equation}
a_{k,1}{x_1}+....+a_{k,r}{x_r}
 \ge   a_{k,1}+...+a_{k_r}-c_k-1
\label{polytope4}
\end{equation}
 Let $(j_1,..,j_r \vert m_1,..,m_r)$
be such that the corresponding vector (\ref{vector})
belongs to the boundary of this polytope. 
Then $1 \in {\cal A}''(j_1,..,j_r \vert m_1,..,m_r)$ and hence
the ideal \newline 
${\cal A}''(j_1,..,j_r \vert m_1,..,m_r)$
is the local ring of the origin (i.e is not proper). 
\par If $\pi: V \rightarrow {\bf C}^2$ is an embedded resolution
then the discrepancy of $f_1^{\gamma_1} \cdot \cdot \cdot f_r^{\gamma_r}$
 along $E_k$ is:
  $$c_k-(a_{k,1}\gamma_1+...+a_{k,r} \gamma_r)$$
i.e. the discrepancy along each $E_k$ is not less than $-1$
if and only if $(1-\gamma_1,...,1-\gamma_r)$ satisfies
(\ref{polytope4}). 

As an example to this proposition one can consider the ordinary cusp $x^2+y^3$
the log-canonical threshold is ${5 \over 6}$ and the constant 
of quasiadjunction is ${1 \over 6}$ (cf. section \ref{alexpolynomialposition}).

The polytopes of quasiadjunction are ``pieces'' of the zeros 
of the (multivariable) Alexander polynomail and in this sense 
are analogs of the spectrum. 

\begin{problem}
Find a generalization of the semicontinuity of spectrum 
of a singularity
\end{problem} 

For some results in this direction cf. \cite{alexhodge}.

\subsection{Hodge decomposition for INNCs.}

The calculation from the last section can be partially extended to 
INNC. Namely we extend the calculation which involve the forms of 
top degree and hence will obtain at least a part of the components
of characteristic varieties. 

We shall start with the definition:

\begin{dfn}(cf. \cite{alexhodge}) Let $f_i=0$ be the equation of divisor $D_i$
and let $\pi: \tilde {\bf C}^{n+1} \rightarrow {\bf C}^{n+1}$ 
be a resolution of the singularities of $\bigcup D_i$ 
(i.e. the proper preimage of the latter is a normal crossings divisor).
Let ${\cal V}_{m_1,...,m_r}$ be the singularity (\ref{icis}) 
having $V_{m_1,...,m_r}$ as its link.  
Let $\tilde V$ be a normalization of 
${\tilde {\bf C}^{n+1}} \times_{{\bf C}^{n+1}} \mathcal{V}_{m_1,..,m_r}$ 
(cf. (\ref{icis}))
The ideal of quasiadjunction of type $(j_1,...,j_r \vert m_1,...,m_r)$
is the ideal ${\cal A}(j_1,...,j_r \vert m_1,....,m_r)$
of germs $\phi \in \mathcal{O}_{0,{\bf C}^{n+1}}$ such that the $(n+1)$-form:
\begin{equation}\label{diffform}
\omega_{\phi}={{\phi z_1^{j_1} \cdot ...\cdot z_r^{j_r} dx_1 \wedge ...
\wedge dx_{n+1}} \over {z_1^{m_1-1} \cdot ... \cdot z_r^{m_r-1}}}
\end{equation}
on the non singular locus of $V_{m_1,...,m_r}$  after the pull back on 
$\tilde V$ extends over the exceptional set.

The $l$-th ideal of log-quasiadjunction  
${\cal A}_l({\rm log} E)(j_1,...,j_r \vert m_1,....,m_r)$
is the ideal of $\phi \in {\cal O}_{0,{\bf C}^{n+1}}$  such that 
the the pull back of the corresponding form $\omega_{\phi}$ 
on $\tilde V$ is ${\rm log}$-form on $(\tilde V, E)$ 
having weight at most $l$.
\end{dfn}

We have the following:

\begin{prop}\label{quasiadjunction}(cf. \cite{alexhodge}, \cite{innc})
There exist a collection of subsets 
${\cal P}_{\kappa}, (\kappa \in \mathcal{K})$
in the unit cube $${\cal U}=\{(x_1,....,x_r) \vert 0 \le x_i \le 1 \}$$
in ${\bf R}^r$ and a collection of affine hyperplanes 
${l}_i(x_1,...,x_r)=\alpha_i$ such that 
each ${\cal P}_{\kappa}$ 
is the boundary of the polytope consisting of solutions to
 the system of inequalities:
$${ l}_i(x_1,...,x_r) \ge \alpha_i$$ 
and such that 
\begin{equation}
({{j_1+1} \over {m_1}}, ... ,{{j_r+1} \over {m_r}}) \in {\cal U}
\end{equation} 
belongs to ${\cal P}_{\kappa}$ if and only if
 
\begin{equation}\label{loginequality} 
 {\rm dim} {\cal A}({\rm log} E)(j_1,...,j_r \vert m_1,....,m_r)/
{\cal A}(j_1,...,j_r \vert m_1,....,m_r) \ge 1 
\end{equation} 

Moreover 
\begin{equation}\label{kloginequality} 
 {\rm dim} {\cal A}_{l}({\rm log} E)(j_1,...,j_r \vert m_1,....,m_r)/
{\cal A}_{l-1}(j_1,...,j_r \vert m_1,....,m_r) \ge k 
\end{equation}
if only if (\ref{vector}) belongs to a collection of certain faces 
${\cal P}_{\kappa,\iota}^{k,l} (\iota \in {\cal I}^{k,l})$
of polytopes ${\cal P}_{\kappa}$.
\end{prop}

Now the exponents of the polytopes of quasiadjunction 
land in the characteristic varieties. More precisely
we have:

\begin{theo}(cf. \cite{innc})
 A character of 
$\pi_1(\partial B_{\epsilon}-\partial B_{\epsilon} \cap 
(\bigcup_{1 \le i \le r} D_i)$ acting on 
\newline $W_l(F^{n}H^n(V_{m_1,....,m_r}))$ via the action of the 
Galois group has the eigenspace of dimension at least $k$
if and only if it has the form: 
$$(exp 2 \pi \sqrt {-1} a_1,...,exp 2 \pi \sqrt {-1} a_r)$$
where $(a_1,....,a_r)$ belongs to one of the 
faces ${\cal P}^{k,l}_{\kappa,\iota}$ of a polytope ${\cal P}_{\kappa}$ 
of quasiadjunction of $\bigcup D_i$. In particular, 
the Zariski closures of exponents of polytopes of quasiadjunction 
are components of characteristic varieties. These components are
the translated subgroups by points of finite order
\end{theo}

We conjecture that all components are the translated subgroups by 
points of finite order. 

\begin{conj}\label{localconjecture}
 Characteristic variety is a union of translated subtori 
of \newline
${\rm Spec}{\bf C}[\pi_1(\partial B_{\epsilon}-\partial B_{\epsilon} \cap X)]$ 
with each translations given by a point of finite order. 
\end{conj}

An interesting problem is to calculate 
them in terms of resolution.

\section{Homotopy groups of the complements to 
hypersurfaces in projective space and linear systems determined 
by singularities}\label{globalalexander}

In this section we want to discuss the characteristic 
varieties associated with hypersurfaces which are divisors
with isolated non normal crossings in a projective space. 
An interesting case occurs already  
when all hypersurfaces have degree one i.e. the case of
arrangements of hyperplanes. The advantage of the case of INNC 
is that one does not have the problems associated with 
complexity of the fundamental group since the fundamental groups for 
such arrangements are abelian (unless we are dealing with an arrangement 
of lines). The theory of such arrangement is still highly non trivial 
and is far from being well understood.
Note that a more general case of divisors with normal crossings 
in general projective manifolds (rather then in $\PP^{n+1}$)
is considered in \cite{sapporo}.
The main results and conjectures of this section show how 
the local characteristic varieties plus certain linear 
systems associated with the points of non normal crossings
determine the global characteristic varieties. This generalizes
the results on the Alexander polynomial discussed earlier.

\subsection{Homotopy groups of the complements to 
INNC}

Let us consider the a divisor $D$ in projective space $\PP^{n+1}$ which 
is a divisor with isolated non normal crossings. This situation 
includes as its special cases:

\bigskip

a) Arbitrary reduced curves in $\PP^2$

\smallskip
b) Hypersurfaces in $\PP^{n+1}$ with isolated singularities and 
hypersurfaces in $\CC^{n+1}$ with isolated singularities and 
transversal to the hyperplanes at infinity.

\smallskip
c) Arrangements of hyperplanes in $\PP^{n+1}$ such that each intersection 
of hyperplanes having codimension $k \ne n+1$ 
belongs to exactly $k$ hyperplanes of the arrangement.  

\bigskip

The starting point is the the following vanishing of 
the homotopy groups generalizing already discussed result from \cite{annals}:

\begin{theo}\label{arrangements}(cf. \cite{sapporo}) 
  Let $X={\bf P}^{n+1}$ and $D$ be a divisor 
having finitely many non-normal crossings.
Assume  that  one of the components has degree 1. 
Then $\pi_i({\bf P}^{n+1}-D)=0$ for $2 \le i \le n-1$. 
If all intersections are the normal crossings, then  
$\pi_n(\PP^{n+1}-D)=0$ and hence $\PP^{n+1}-D$ is 
homotopy equivalent to the wedge of the $n+1$-skeleton of 
the torus $(S^1)^k$ and several copies of $S^{n+1}$. 
\end{theo}

One also has a similar vanishing for the homology of local systems.

\begin{theo} \label{locsys}(cf. \cite{sapporo})
 Let $\chi \in {\rm Char}\pi_1({\bf P}^{n+1}-D)$ be a 
character of the fundamental group different from the identity and 
let ${\CC}_{\chi}$ be $\CC$ considered as $\CC [\pi_1({\bf P}^{n+1}-D)]$ 
module via the character $\chi$. Then 
$$H_i({\bf P}^{n+1}-D,\chi)=0 \ (i \le n-1)$$  
$$H_n({\bf P}^{n+1}-D,\chi)=\pi_n({\bf P}^{n+1}-D) 
\otimes_{\CC [\pi_1({\bf P}^{n+1}-D)]} {\CC}_{\chi}$$
\end{theo} 

The main problem for INNC hence is to understand 
the first non trivial homotopy group
$\pi_n(\PP^{n+1}-D)$. Similarly to the local case the starting 
point is the following:

\begin{dfn}(cf. \cite{sapporo})
 The $k$-th characteristic variety $V_k(\pi_n({\bf P}^{n+1}-D))$
of the homotopy group $\pi_n({\bf P}^{n+1}-D)$
is the zero set of the $k$-th Fitting ideal of $\pi_n({\bf P}^{n+1}-D)$, i.e. 
the zero set of minors of order $(n-k+1) \times (n-k+1)$ of $\Phi$ in a 
presentation 
$$\Phi: \CC [\pi_1({\bf P}^{n+1}-D)]^m \rightarrow
 \CC [\pi_1({\bf P}^{n+1}-D)]^l 
\rightarrow \pi_n(X) \rightarrow 0$$
 of $\pi_1({\bf P}^{n+1}-D)$ module $\pi_n({\bf P}^{n+1}-D)$ 
via generators and relations. 
Alternatively (cf. theorem \ref{locsys})
outside of $\chi=1$, $V_k(\pi_n({\bf P}^{n+1}-D))$ is the set 
of characters $\chi \in \Char [\pi_1({\bf P}^{n+1}-D)]$ such that  
$\dim H_n({\bf P}^{n+1}-D,\chi) \ge k$.
\end{dfn}

\subsection{Jumping loci on quasiprojective varieties} 

A remarkable fact is that the characteristic varieties of 
the complements have a very simple structure (unlike in 
the similar situations outside of algebraic geometry). We did see this already 
in the case of links of curve singularities and in the case of charactersstic 
variety $V_1$ in general local case. 
The local systems on {\it a non singular projective 
variety} correspond to holomorphic bundles which are topologically 
trivial. The jumping loci for the cohomology of such bundles are unions 
of translates of abelian subvarieties of the Picard variety. These 
results, having long history, are due to Catanese, Beauville, 
Green-Lazarsfeld, Simpson and Deligne. 
We shall use the following quasi-projective version dealing with the
cohomology of local systems which will allow us eventually to 
describe the characteristic varieties.

\begin{theo} (\cite{arapura}) Let $\hat X$
 be a projective manifold such that $H^1(\hat X,\CC)=0$.
Let $\hat D$ be a divisor with normal crossings. Then there
exists a finite number of unitary characters $\rho_j \in {\rm Char} 
\pi_1(\hat X-\hat D)$ and  holomorphic maps $f_j: \hat X-\hat D
\rightarrow T_j$ into complex tori $T_j$
such that the set 
$\Sigma^k(\hat X-\hat D)=\{ \rho \in {\rm Char} \pi_1(\hat X-\hat D)
\vert {\rm dim} H^k(\hat X-\hat D,\rho)  \ge 1  \} $
coincides with $\bigcup \rho_jf^*_iH^1(T_j,\CC^*)$. In particular,
$\Sigma^k$ is a union of translated by unitary characters 
subgroups of ${\rm Char}\pi_1(\hat X-\hat D)$.
\end{theo} 
Hence we also obtain:

\begin{coro}\label{arapurahomotopy } 
The characteristic variety $V_k(\pi_n(\PP^{n+1}-D))$ 
is a union of translated subgroups $S_j$ of the group 
${\rm Char} \pi_1(\PP^{n+1}-D)$
by unitary characters $\rho_j$: $$V_k(\pi_n(\PP^{n+1}-D))=\bigcup \rho_jS_j$$
\end{coro}

In the case $k=1$ the components of characteristic variety having 
a positive dimension correspond to the maps onto hyperbolic 
curves. This has many applications for example to calculations 
of characteristic varieties (cf. \cite{charvar}),
 estimating the order of the group
of automorphisms of the complements (cf. \cite{bandman}), 
classification of arrangements of lines (cf.\cite{sergey}) among others
but we won't discuss them here.

\subsection{The Hodge numbers of abelian covers 
of projective spaces and 
linear systems.} 
Let, as before, 
$D=\cup_{i=0,..r} D_i$ be a divisor in $\PP^{n+1}$.
We shall assume, to simplify the exposition, that one of components,
say $D_0$ has degree equal to one  and that there are $D$ has non
non normal crossings on $D_0$.
Let $\pi_1(\PP^{n+1}-D) \rightarrow \oplus \ZZ/m_i\ZZ$ 
be a surjective homomorphism and
let $X_{\bf m} \ ({\bf m}=(...,m_i,...))$ 
be a normalization of a compactification of 
unbranched cover of $\PP^{n+1}-D$ corresponding to this homomorphism.
Let $f: X_{\bf m} \rightarrow \PP^{n+1}$ be the corresponding projection.

Our goal is to calculate the Hodge number $h^{n,0}(X_{\bf m})$.
Starting from $D$, we shall define global polytopes 
of quasiadjucntion so that with each face $\delta$ of the 
polytope is associated the ideal sheaf ${\cal J}_{\delta}$. 
The above Hodge number is equal 
to the number of lattice points is $\delta$ 
counted with the weight given 
by the dimension of linear system of hypersurface of degree given by 
$\delta$ and with local conditions given by the ideal sheaf 
${\cal J}_{\delta}$.

Let is consider the unit cube $\U=\{ (x_1,...,x_r) \in \RR^r \vert 
0 \le x_i \le 1\}$ coordinate of which correspond to  
irreducible components $D_1,...,D_r$ of the divisor $D$.
We view $\RR^r$ as the universal cover of the group $(S^1)^r$
of unitary characters of $\pi_1(\PP^{n+1}-D)$ and 
$\U$ as the fundamental domain for the action of the 
covering group on the cover.
With each point $P$ in $\PP^{n+1}$ where $D$ has 
a non-normal crossing the definition \ref{quasiadjunction}
associates a polytope $\P_{\kappa}$
in the unit cube in $\RR^s$ with coordinates corresponding 
to the components of $D$. Since one has the canonical projection
$\pi: \RR^r \rightarrow \RR^s$, forgetting the coordinates 
corresponding to $D_i$'s not containing $P$, each $\P_{\kappa}$ 
defines the polytope $\pi^{-1}(\P_{\kappa})$ in $\U$ which 
we shall denote by the same letter. This defines a finite 
collection of polytopes $\P_{\kappa,P} \subset \U$.

\begin{dfn}\label{globalquasiadjunction} 
Consider the equivalence relation on points  in $\U$ calling two 
points equivalent if the collections of polytopes $\P_{\kappa,P}$ 
containing these two points are identical. 
The equivalence class is called the global polytope of quasiadjucntion.

\smallskip

\noindent 
A global face of quasiadjunction is a face of a global polytope of 
quasiadjunction. 

\smallskip

\noindent Let $\S_{\delta}$ be the set of non normal crossings $P$ of $D$ for 
which there exist the polytopes $\P_{\kappa,P}$
 contaning $\delta$. 

\smallskip
\noindent 
The ideal sheaf corresponding to $\delta$ is a sheaf 
$\J_{\delta} \subset \O_{\PP^{n+1}}$ such that 
$\O_{\PP^{n+1}}/\J_{delta}$ is supported at $\S_{\delta}$ and such that 
the stalk 
at $P$ is the ideal which is the intersection of local
ideals of quasiadjunction corresponding to local polytopes containing $\delta$.
\end{dfn}

Clearly such an equivalence class is a polytope i.e. consists
of points satisfying a system of linear inequalities. 
Also the collections $S_{\delta}$ of non normal crossings 
are defined entierly by the local data of $D$.
The Hodge number $h^{n,0}(X_{\bf m})$ depends on 
additional piece of information.

\begin{theo}\label{branchedcoverslinearsystems}(cf. \cite{charvar},
 \cite{sapporo})
Let $D$ as above and let $d_i$ be the degree of the irreducible component
$D_i$.
For each $\chi \in \Char \oplus \ZZ/m_i\ZZ \subset \Char\pi_1(\PP^{n+1}-D)$ 
let $\delta(\chi)$ be the global face of quasiadjunction containing 
${1 \over {2 \pi \sqrt{-1}}}{\rm log}(\chi) \in \U$. Let $l$ be such that 
the hyperplane $d_1x_1+...+d_rx_r=l \ (l \in \ZZ)$ contains $\delta$. Then 
$$h^{n,0}(X_{\bf m})=\sum_{\chi \in \Char \oplus \ZZ/m_i\ZZ} 
{\rm dim}H^1(\PP^{n+1},\J_{delta}(l-n-1)$$
\end{theo}

A proof in cyclic case and in the case of curves and 
generalizing Zariski's approach 
(\cite{irreg}) is given in \cite{position} and \cite{charvar}
and the case of INNC is similar.
Alternatively, one can also use the approach in \cite{esnault}.

\begin{exam}
For an irreducible curve of degree
$d$ with nodes and the ordinary cusps as the only singularities 
the global polytope of quasiadjunction coincide with the local one of the 
cusp. The only face of quasiadjunction is $x={1 \over 6}$. The contributing
hyperplane is given by $dx={d \over 6}$ and its level is $d \over 6$.
The sheaf of quasiadjunction corresponding to this face of quasiadjunction
is the ideal sheaf having stalks different from the local ring only 
at the points of ${\bf P}^2$ where the curve has cusps and the 
stalks at those points are the maximal ideals of the 
corresponding local rings.
\end{exam}

For characters not on the global faces of quasiadjunction one still 
can define the ideal sheaves looking at the polytopes containing 
the lifts of the characters ${1 \over {2 \pi \sqrt{-1}}}{\rm log}(\chi) \in \U$
into universal cover of the torus of unitary character 
and also the integer $l$ such that $d_1x_1+...+d_rx_r=l$ 
contains the lift. However the 
corresopnding group $H^1(\PP^{n+1},\J_{delta}(l-n-1)$
will be vanishing. 
For plane curves with nodes and cusps one obtains the 
following classical result (for the most part already discussed earlier).

\begin{coro}(Zariski's theorem) Let $C$ be a plane curve of degree $d$ 
having nodes and cusps as the only singularities. Let $\J$ be a subsheaf of 
the sheaf of regular functions whose sections belong to the maximal 
ideals at the points in $\PP^2$ which are cusps of $C$.
If $k > {{5d} \over 6}$
then $$H^1(\PP^2, \J(k-3))=0$$
If $6 \vert d$ then $H^1(\PP^2, \J({{5d} \over 6}-3))=h^{1,0}(X_d)$ 
is equal to the 
irregularity of a a resolution of singularities of a 
$d$-fold cyclic cover of $\PP^2$ branched over $C$.
\end{coro}

\subsection{Mixed Hodge structure on homotopy groups}

The theorem \ref{locsys} suggests an additional structure on  
the characteristic varieties coming from the mixed Hodge structure
on the cohomology of local systems. This is an analog of 
discussed earlier in local case the Hodge decomposition of 
characterstic varities.
The MHS on the cohomology of local systems can be understood 
by interpreting the cohomology of local systems having finite order 
as the eigenspaces of the Galois group acting on the abelian covers 
as follows.

\begin{theo}\label{hodgelocalsystems}
Let $G$ be a finite group and $g: \pi_1(X) \rightarrow G$ be a surjection.
Let $\chi \in \Char \pi_1(X)$ which is the pull back of a character 
of $G$. Assume that $\pi_i(X)=0$ for $2 \le i \le n-1$. Finally 
let $X_G$ be the unbranched cover of $X$ corresponding to $g$. 
Then the eigenspace $H^n(X_G)_{\chi}$ is isomorphic to the homology of 
$H^n(\CC_{\chi})$ of the local system $\CC_{\chi}$ corresponding to $\chi$.
In particualar, the cohomology classes in $H^n(\CC_{\chi})$ acquire
 the Hodge type. 
\end{theo}

If $X$ is quasiprojective and non singular, so is $X_G$ and hence $H^n(X_G)$ 
admits the mixed Hodge structure with the weights $n,....,2n$. 

\begin{dfn}(cf. \cite{sapporo}) 
Let $\PP^{n+1}-D$ be a complement to an INNC in $\PP^{n+1}$.
For a local system $\chi$ pf finite order 
let $h^{p,q,n})_{\chi}$ be the dimension of the space of cohomology classes
in $H^n(\CC_{\chi})$ having the Hodge type $(p,q)$.
The following subset of $V_k(\pi_n(\PP^{n+1}-D))$:
$$\P^{p,q,n}_k=\{ \chi \vert h^{p,q,n}_{\chi} \ge k \}$$
is called the component of the characteristic variety of type 
$(p,q,n)$
\end{dfn}

One has $\P^{p,q,n}_k \ne \emptyset$ only if $n \le p+q \le 2n$ 
and $\bigcup_{p,q,n}\P^l=V_1(\pi_1(\PP^{n+1}-D))$.

\subsection{A relation between the Hodge numbers of 
branched and unbranched abelian covers}

We want to use the theorem \ref{branchedcoverslinearsystems}
to detect some components of characteristic varieties of the homotopy 
groups. Here is the relation between branched and unbranched covers
which we shall need since
the theorem \ref{branchedcoverslinearsystems} works in compact case.

\begin{theo}\label{branchedcover}(cf. \cite{sapporo})
Let
$\chi \in Char(\pi_1(\PP^{n+1}-D))$ be a character of 
a finite quotient $G$ of $\pi_1(\PP^{n+1}-D)$. Let $\bar U_G$ 
be a $G$-equivariant non-singular compactification of $U_G$ and
let $H^{p,q}(\bar U_G)_{\chi}$ be the $\chi$-eigenspace of $G$ 
acting on $H^{p,q}(\bar U_G)$. 
Then 
$$ h^{n,0,n}(\CC_{\chi})=h^{n,0}(U_G)_{\chi}=
h^{n,0}(\bar U_G)_{\chi}$$
\end{theo}

\subsection{Main theorem and Open Problems}

Combining this together we obtain the following, extending 
results on Alexander polynomials and the case of reducible 
curves in \cite{charvar}:

\begin{theo}\label{charvarposition}(cf. \cite{sapporo})
 Let $D \subset \PP^{n+1}$ be a union of hypersurfaces
$D_0,D_1,...,D_r$ of degrees $1,d_1,...,d_r$ respectively, 
which is a divisor with isolated 
non-normal crossings. Let $\F$ be a face of global polytope of 
quasi-adjunction, i.e. 
a face of an intersection of polytopes of quasi-adjunction corresponding 
to a collection $\cal S$ of 
 non-normal crossings of $D$. Let $d_1x_1+...+d_rx_r=l$ be
a hyperplane containing the face of quasiadjunction $\F$.
If $H^1(\A_{\F}\otimes \O(l-3)) =k$, then the Zariski closure of 
$exp (\F) \subset \Char H_1(\PP^{n+1}-D)$ belongs to a component of 
$V_k(\pi_n(\PP^{n+1}-D))$.
 
\end{theo}

There is a generalization to INNC divisors on arbitrary 
projective simply-connected varieties. I refer to \cite{sapporo}
for conjectures. Here is a short list of the open problems
even in the case of divisors in $\PP^{n+1}$. 

\begin{problem} Are there components of characteristic variety 
$V_k(\pi_1(\PP^{n+1}-D))$ which are not Zariski closures of 
$\P_{n,0,n}$?
\end{problem}

\begin{problem} Find methods for detecting the sets 
$\P_{p,q,n}$ with $(p,q) \ne (n,0)$
\end{problem}

A difficulty here is that one cannot work with arbitrary 
compactification since the Hodge numbers $h^{p,q}$ are not 
birational invariants. It would be good to have techniques 
which will allow to work directly with the complement and avoiding
to some extent the compactification.

\begin{problem} Generalize the main theorem
to projective algebraic varieties and  beyond the cases when 
$\O(D_i)=\L^{m_i}$.
\end{problem}

See discussion of this in \cite{sapporo}

\begin{problem} Find additional interesting examples
beyond the one described in \cite{sapporo}.
\end{problem} 

\begin{problem} Solve realization problem for characteristic varieties
i.e. describe how many components and what are their dimensions
depending on the numerical data of the divisor $D$ on $X$.
\end{problem}




\end{document}